\documentclass[aps,twocolumn,floatfix,nofootinbib]{revtex4-1}% standard PRL
\usepackage{dcolumn}% Align table columns on decimal point
\usepackage{times}
\usepackage{graphics}
\usepackage{bm}
\usepackage{subfigure}
\usepackage[pdftex]{graphicx}
\usepackage{amsthm}
\usepackage{notes2bib}
\usepackage{amssymb}
\usepackage{url}
\usepackage{enumerate}
\usepackage{epstopdf}
\usepackage{mathtools}          
\usepackage{setspace}
\usepackage[usenames,dvipsnames,svgnames,table]{xcolor}
\usepackage{color}
\usepackage{algorithmic}
\usepackage{algorithm}

\newcommand{\allblack}{\color{black}{}}

\begin{document}

\preprint{APS/123-QED}

\title{
%A modern approach of Isaac Newton:\\
Constructing low-dimensional ordinary differential equations from chaotic time series\\ of high\slash infinite-dimensional systems using radial function-based regression}% Force line breaks with \\

\author{Natsuki Tsutsumi}
\affiliation{Faculty of Commerce and Management, Hitotsubashi University, Japan}
\author{Kengo Nakai}
\affiliation{The Graduate School of Environment, Life, Natural Science and Technology, Okayama University, Japan}
\author{Yoshitaka Saiki}
\affiliation{Graduate School of Business Administration, Hitotsubashi University, Japan}
%\keywords{data-driven modeling $|$ linear regression $|$  ordinary differential equations $|$dynamical systems $|$ chaos}
\date{\today}% It is always \today, today,
\begin{abstract}
In our previous study (N. Tsutsumi, K. Nakai and Y. Saiki (2022)) we proposed a method of constructing a system of differential equations of chaotic behavior only from observable deterministic time series, which we will call radial function-based regression (RfR) method. 
The RfR method employs a regression using Gaussian radial basis functions together with polynomial terms to facilitate the robust modeling of chaotic behavior.
In this paper, we apply the RfR method to several types of relatively high-dimensional  deterministic time series generated by a partial differential equation,  a delay differential equation, a turbulence model, and intermittent dynamics.
The case when the observation includes noise is also tested.
We have effectively constructed a system of differential equations for each of these examples,  which is assessed from the point of view of time series forecast, reconstruction of invariant sets, and invariant densities. 
{
We find that in some of the models, an appropriate trajectory is realized on the chaotic saddle and is identified by the Stagger-and-Step method.}
\end{abstract}

\maketitle
%\thispagestyle{firststyle}
%\ifthenelse{\boolean{shortarticle}}{\ifthenelse{\boolean{singlecolumn}}{\abscontentformatted}{\abscontent}}{}

%%%%%%%%%%%%%%%%%%%%%%%%%%%%%%%%%%%%%%%%%%%%%%%%%%%%%%%%
\section{Introduction}
Scientists have attempted to search for a governing law  
when they observe an 
intriguing phenomenon.
Kepler’s laws of planetary motion 
are said to be derived by 
Johannes Kepler, whose analysis of the observations of 
Tycho Brahe enabled him to establish the laws in the early 17th century.
Sir Isaac Newton is said to establish the law of universal gravitation based on experimental observations made previously by Galileo Galilei.
The laws are summarized as differential equations. 
Since then especially for microscopic dynamics of various physical phenomena, the corresponding differential equations
have been discovered.

Models of dynamics are constructed 
primarily via a physical understanding of the phenomena. 
However, in the last few decades, several approaches have been proposed concerning modeling dynamics from given time series data with the aid of machine learning techniques~\cite{chen2018,Pathak_2017,nakai_2018,nakai_2020,kobayashi2021}.
Although the aim of these approaches is mainly to infer short time series, some models have also succeeded in mimicking dynamical system features 
such as invariant sets and invariant densities~\cite{kobayashi2021}.

Some studies~\cite{baake1992,wang2011,brunton16, champion19}
estimate a system of ordinary differential equations (ODEs) {from time series data}, which makes it easy to analyze dynamical system features.
The above methods require observation time series of all the variables of a system of ODEs to be modeled. 
For example, they require time series of all three variables $x, y, z$ to construct a data-driven model of the Lorenz 1963 dynamics~\cite{lorenz_1963}.

There are attempts~\cite{gouesbet91,gouesbet92,gouesbet97} which derive a system of ODEs from scalar time series. 
The observable variable and its time derivatives are used for the model variables.
 In \cite{petrov03}, 
the problem 
with observational (Gaussian or white) noise is investigated.
In the series of studies, prior knowledge of the background dynamics is used to choose basis functions for the regression, and the method is not appropriate for modeling for practical purposes.

Recently we proposed a simple method of constructing a system of ODEs of chaotic behavior based on the regression only from observable scalar time series data with the basis functions unchanged
~\cite{tsutsumi22}.
Independent of the dynamics, we employ 
 spatially localized radial basis functions
in addition to polynomial basis functions for the regression. 
We can say that 
our method is the only method that can model a system of ODEs from scalar time series data which is applicable to various dynamics 
without assuming the knowledge of the background dynamics.
Furthermore, we should note that 
the capability of the forecasting is not limited to a short time series but a density distribution created using a long time series.

Unlike
prior studies~\cite{baake1992, wang2011, brunton16} 
our radial function-based regression (RfR) method proposed in~\cite{tsutsumi22} is applicable 
even when the time derivative of each variable 
is not approximated by the low-order polynomials of variables. 
We exemplified that our RfR method worked well for 
a time series of one variable ($x$) of the typical and simple chaotic dynamics given by the Lorenz system, and evaluated the constructed data-driven model in detail.
Our model construction using the RfR method has several advantages in that the method simultaneously meets the following: 
a model variable is physically understandable; 
a model can be constructed even when the number of observable variables is limited
and even when no knowledge of the governing system is given. 
The RfR method allows us to construct a model 
using physically understandable variables.

In this research, we present the applicability of our RfR method to complex dynamics described by infinite or finite but high-dimensional systems including 
the Kuramoto-Sivashinsky equation~\cite{christiansen97},
the Mackey-Glass equation~\cite{mackey77}, the shell model of fluid turbulence~\cite{yamada87}
and intermittent dynamics described by the coupled R\"{o}ssler equation~\cite{pikovsky02}.
The Kuramoto-Sivashinsky equation is the partial differential equation system that can reveal spatio-temporal chaos; the Mackey-Glass equation is the delay differential equation which is a model of feedback control of blood cells; 
the shell model turbulence is a system of
high-dimensional ODEs of complex variables which mimic the spectral equation in Fourier space of the Navier-Stokes system; the coupled R\"ossler equation is 
an ordinary differential equation system that can show intermittent behavior.    
Note that for the shell model turbulence, 
we use time series data of a variable that does not appear in the original system.
It is common that observable data includes noise. 
So in the case of the Kuramoto-Sivashinsky equation, the effect of noise added to the observation data is also investigated.

We assume there exists an  unknown 
system of $N$ dimensional 
ODEs called an original system concerning an unknown variable $\bm{x}$:
\begin{equation}
 \frac{d\bm{x}}{dt}=\bm{f}(\bm{x}),
 \label{eq:odegeneral}
\end{equation}
or an infinite dimensional system such as a partial differential equation and delay differential equation  whose dynamics can be approximated well by Eq.~\ref{eq:odegeneral}. 
In this paper the above-mentioned examples are chosen as Eq.~\ref{eq:odegeneral}.
We can observe some of the components of the variable $\bm{x}$,
or more generally 
\begin{equation}
    w_i=g_i(\bm{x}),~ i=1,\ldots, I. \label{eq:observable} 
\end{equation}
Note that $w_i$ does not necessarily be described as a function of $\bm{x}$, 
but can be an integrated value of $\bm{x}.$
The case when $w_i$ includes noise is also investigated.
We are unable to reconstruct the original system Eq.~\ref{eq:odegeneral} itself from given time series data unless all the components of the variable $\bm{x}$ in Eq.~\ref{eq:odegeneral} and their time series are known. 
We assume that there exists 
a system of $D$ dimensional 
ODEs:
\begin{equation}
    \frac{d\bm{X}}{dt} = \bm{F}(\bm{X}),
    \label{eq:odeestimation}  
\end{equation}
where the first $I$ components of the variable $\bm{X}$ are $X_i=w_i$ $(i=1,~\ldots,I)$ and the rest of $X_i$ $(i=I+1,~\ldots,D)$ are created from the delay-coordinates of some $w_i$~\cite{takens_1981,sauer1991}.
The aim of our proposed method is to model~Eq.~\ref{eq:odeestimation}
that can describe the behavior of the observable variables $w_i$ $(i=1,~\ldots,I)$ 
as components of the variable $\bm{X}$ 
which is made only by the observable variable $w_i$.

By the RfR method, we
have succeeded in the modeling of 
dynamics 
even when the right-hand side of Eq.~\ref{eq:odeestimation} cannot be described by low-order polynomials.
This is primarily because of the introduction of the Gaussian radial basis functions.
It allows us to ease the modeling of a chaotic attractor for a variety of dynamics.
However, even if the obtained data-driven model does not have a chaotic attractor, we can generate a trajectory staying in the chaotic saddle 
by employing the Stagger-and-Step method~\cite{sweet_2001c}.
The Stagger-and-Step method is applied to create a long chaotic trajectory of the data-driven model for the shell model turbulence and the Kuramoto-Sivashinsky equation with observation noise. 
See Appendix~\ref{sec:Stagger-and-Step} for the detail.

The rest of the paper is organized as follows. 
In Section~\ref{sec:method}, we introduce the proposed method for deriving a system of ODEs.
In Section~\ref{sec:model}, we demonstrate the settings of the construction of a data-driven model for each of the examples.
In Section~\ref{sec:results}, we assess the constructed data-driven models. 
Concluding remarks are given in Section~\ref{sec:conclusion}.
\section{Method: Radial function-based regression (RfR)}\label{sec:method}
\begin{figure}%[ht]
\vspace*{2mm}
  \begin{center}
    \includegraphics[width=0.99\columnwidth,height=0.65\columnwidth]{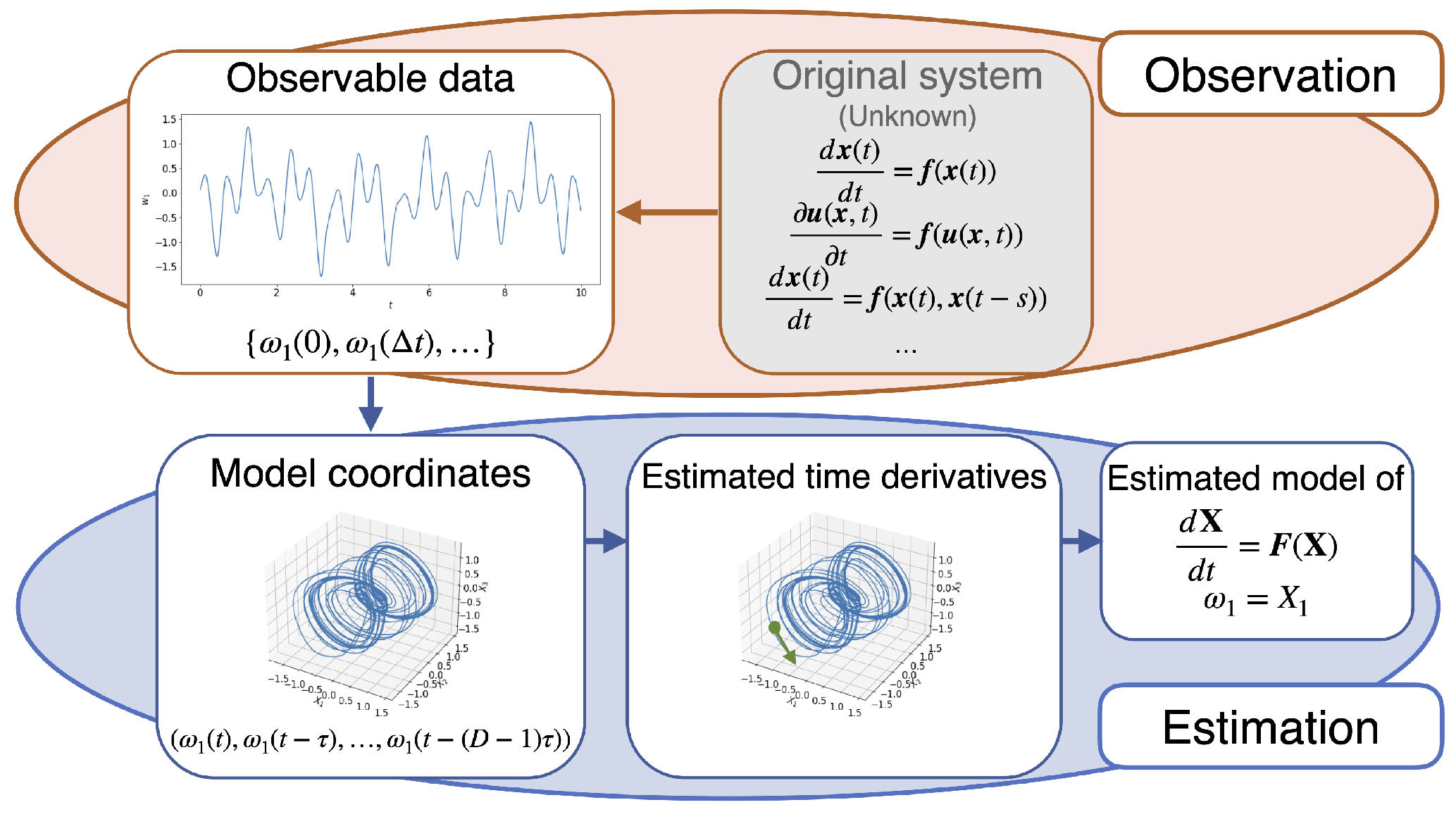}  
   \end{center}
 \caption{
{\bf Outline of the proposed radial function-based regression (RfR) method for constructing a system of ODEs.} 
 }
 \label{fig:outline}
\end{figure}
As shown in Fig.~\ref{fig:outline}, our aim is to construct a system of $D$ dimensional ODEs~Eq.~\ref{eq:odeestimation} based only on some observable 
deterministic time series $\omega_i (i=1,\ldots ,I)$
of length $T$ with time step $\Delta t$,
and data size $N_T = T / \Delta t$.
The steps of the proposed method 
(RfR)
are outlined below:
\begin{itemize}
    \item[{i.}] Choose the delay-coordinate variable: 
    delay-time $\tau$ and dimension $D$ 
    \item[{ii.}] Estimate the time derivative at each sample point using the Taylor approximation 
    \item[{iii.}] Choose the basis function\footnote{We do not assume some particular vector space in advance as is usual in the literature for machine learning.} used in Step iv 
    \item[{iv.}] Perform linear regression at sample points with ridge regularization 
    \item[{v.}] Assess the model quality according 
    to the reproducibility of the delay structure in a generated model trajectory
\end{itemize}

\subsection{Choice of model variables}\label{sec:delay}
For practical reasons, the number of observable variables is small, when compared to the dimension of the background dynamical system. 
To describe dynamics,
we generate a higher dimensional variable by introducing 
the delay-coordinates of limited observables 
such as 
$(w_1(t),w_1(t-\tau), w_1(t-2\tau),\ldots, w_1(t-(D-1)\tau))$.

\subsection{Estimating the time derivative}
\label{sec:timederivative}
We apply
the Taylor approximation to estimate the time derivative at each sample point $\bm{X}(\tilde{t})$ based only on discrete time points of a trajectory $\bm{X}(t)$ $(t=\ldots, \tilde{t}-\Delta t,\tilde{t},\tilde{t}+\Delta t,\ldots)$.
In our computations, we employ the sixth order approximation:
\begin{align*}
 &\frac{d\bm{X}(\Tilde{t})}{dt} \approx 
 \frac{1}{60 l\Delta t} \{ \bm{X}(\Tilde{t}+3l\Delta t) - 9\bm{X}(\Tilde{t}+2l\Delta t)\\ ~~~&+45 \bm{X}(\Tilde{t}+l\Delta t) 
- 45\bm{X}(\Tilde{t}-l\Delta t) + 9\bm{X}(\Tilde{t}-2l\Delta t)\\ 
~~~&-\bm{X}(\Tilde{t}-3l\Delta t)\},
\end{align*} 
where $l~(\ge 1)$ is a positive integer to determine the points to estimate the time derivative at $\bm{X}(\Tilde{t})$, that is, points at every $l\Delta t$ time step are used.

When the observable data includes noise, we choose a large value of $l$
to estimate the time derivative, which allows us to avoid high-frequency oscillations. 
In that case, we need to employ the high ({\textit{e}.\textit{g}.}, sixth) order Taylor approximation.
Note that if the observable data does not include noise, we do not need to use such a high order formula.
See Appendix~\ref{sec:relation_taylor_time-step} for the estimation of the time derivative at each sample point when noise is added to the time series.

\subsection{Choice of basis function}
\label{sec:basefunctions}
We construct a model Eq.~\ref{eq:odeestimation} through the linear regression of the following form:
\begin{equation}
F_k(\bm{X}) \approx \tilde\beta_0^k + \sum_{d=1,\cdots,D} \tilde\beta_d^k X_d + \sum_{j=1,\cdots,J} \tilde\beta_{D+j}^k~\phi_j(\bm{X}),\label{eq:gaussianpoly}
\end{equation}
where $F_k(\bm{X})$ is the $k$th component of $\bm{F}(\bm{X})$ in Eq.~\ref{eq:odeestimation}, 
${\tilde{\boldsymbol{\beta}}^{k}}=(\tilde{\beta}^k_0,\tilde{\beta}^k_1,\cdots,\tilde{\beta}^k_{D+J})$ is a set of estimated parameters and  
\begin{equation*}
\phi_j(\bm{X}) = \exp\left(\frac{-||\bm{X}-c_j||^2}{\sigma^2}\right),%\label{eq:gaussian}
\end{equation*}
is the Gaussian radial basis function, 
where $\|\cdot\|$ denotes the $l^2$ norm,
$c_j \in \mathbb{R} ^D$ is the coordinate of the $j$~th center point ($j=1,\ldots,J$), and $\sigma^2$ is the parameter that determines the deviation of $\phi_j$.

In this RfR method,
$c_j$ is distributed as lattice points with grid size $\delta_{grid}$.
For a given integer $m$ explained below, 
we consider $c_j$
such that 
there exists a data point in the $(m-1)\delta_{grid}$-neighborhood.
An increase in $\delta_{grid}$ 
results 
in a decrease in the number of center points and a decrease in the required computational resources.

For a given $\delta_{grid}$, $\sigma^2$ is  %obtained from  
determined as 
\begin{align*}
     \sigma^2 := \frac{((m-1) \delta_{grid})^2}{- \log_{e} p},
\end{align*}
where $m$ is the degree of the corresponding B-spline basis function and $p~(>0)$ is a small value. 
When we set $m=3$, and $p=0.1$, then %$\frac{(m-1)^2 }{- \log_{e} p} \approx 1.7372$.
$(m-1)^2 /(- \log_{e} p) \approx 1.7372$.
See Kawano and Konishi~\cite{kawano2007} for more details.

\subsection{Ridge regression}\label{sec:regression}
For determining  
the coefficients in Eq.~\ref{eq:gaussianpoly} 
we obtain the minimizer of the following function 
\begin{align*}
    L(\bm{b}) = \frac{1}{2n} ||\bm{y} - A \bm{b} ||^2 + \frac{\lambda}{2} ||\bm{b}||^2,
\end{align*}
where $n$ is the size of regression data,
$\bm{y}$ represents the 
standardized
time derivative (see the left-hand side of Eq.~\ref{eq:odeestimation}),
$\lambda$ is a parameter determining the strength of regularization, 
and $A$ is $n \times (D+J)$ matrix whose $i$-th row is 
$(X_1(t_i), \ldots X_D(t_i), \phi_1(\bm{X}(t_i)), \ldots, \phi_J(\bm{X}(t_i)) )$ 
(see Eq.~\ref{eq:gaussianpoly}).
The regularization is used to prevent overfitting primarily due to the introduction of the Gaussian radial basis functions.

The ridge estimator attaining the minimum of $L(\bm{b})$ is written as follows: 
\begin{equation*} 
    ({A^\mathrm{T}} A + n \lambda I )^{-1} {A^\mathrm{T}}  \bm{y}, 
\end{equation*}
where $I$ is the identity matrix and $A^\mathrm{T}$ is the transpose of $A$. 
Due to the limited computational resources,
it is hard to use all data of size $N_T$.
We randomly choose samples of size $n$ among $N_T$ for regression.
\subsection{Evaluation of the model}
\label{sec:delaystructure}
For assessing the model quality we confirm the delay structure in the model trajectory. 
When the number of observable variable is one, the model trajectory should satisfy the following: 
\begin{align}
    \label{eq:delay_structure}
    X_1(t) \approx X_2(t+\tau) \approx \cdots \approx X_D(t+(D-1)\tau).
\end{align}
The hyper-parameters such as a regularization parameter and grid size are selected adequately based on the level of the reconstruction of the delay structure.

\subsection{Choice of parameters}\label{sec:parameters}
Table~\ref{table:data_type} depicts parameter settings
that will be employed in the modeling in Sections~\ref{sec:model} and \ref{sec:results}. 
\begin{table}[t]
\small
    \centering
    \begin{tabular}{lccccccccccc}
        \hline 
                    & KS  & MG &  SM & CR & n-KS\\\hline
        $D$         & 5 & 7 & 6 & 6 & 5\\\hline
        $\tau$      & 0.12 & 0.5 & 18 & 0.4 & 0.12 \\\hline
        $m$         & \multicolumn{5}{c}{3}\\\hline
        $p$         & \multicolumn{5}{c}{0.1}\\\hline
        $\delta_{grid}$& 0.50 & 0.25 & 0.25 & 0.25 & 0.50 \\\hline
        $J$  & 19,322 & 241,402 & 693,749 & 91,393 & 35,179 \\\hline  
        $l$         & 1 & 1 & 1 & 1 & 9\\\hline
        $\lambda$   & $10^{-7}$ & $10^{-7}$ & $10^{-12}$ & $10^{-7}$ & $10^{-4}$ \\\hline
        $n$      & \multicolumn{5}{c}{50,000} \\\hline 
        $I$         & 1 & 1 & 1 & 2 & 1\\\hline
        $N_T$     & \multicolumn{5}{c}{$10^6$}\\\hline
        $\Delta t$  & 0.01 & 0.01 & 1.0 & 0.1 & 0.01\\\hline
    \end{tabular}
    \caption{{\bf Sets of parameters.}
    They are used for the modelings of the Kuramoto-Sivashinsky  equation~(KS), Mackey-Glass equation~(MG), the shell model of fluid turbulence~(SM), the coupled R\"ossler equation~(CR), and the Kuramoto-Sivashinsky equation with noise~(n-KS).
    }
\label{table:data_type}
    \normalsize
\end{table}
The parameter $\delta_{grid}$ should correspond to a scale of the variable $\omega_i$ in Eq.~\ref{eq:observable}.
Hence,
$\omega_i$ is standardized in our modeling to avoid the adjustment.
For the choice of parameters we consider the followings:
\begin{itemize}
    \item The dimension $D$ of the model is selected to be bigger than the expected attractor dimension. But 
    small dimension $D$ is effective for the computations.
    \item The delay time $\tau$ is selected based on the decay of correlation of a variable $X_1(t)$. $\tau$ is chosen to meet that the correlation between $X_1(t)$ and $X_1(t-\tau)$ is around 0.5 (See Appendix~\ref{sec:auto-corr}).
    \item The parameter $l$ used for estimating the time derivative at each sample point is chosen larger when the observable data includes noise with larger amplitude (See { section}~\ref{sec:timederivative}).
    \item The number of center points $J$ is determined by the settings of $\delta_{grid}$ and $m$
    (See section~\ref{sec:basefunctions}).
     \item The number of regression points $n$ is not very sensitive and fixed as 50,000 (See section~\ref{sec:regression}).
    \item The regularization parameter $\lambda$ is chosen so that the reproducibility of the delay structure is successful (See section~\ref{sec:delaystructure}), but is not very sensitive.
\end{itemize}

\section{Example dynamics to be modeled}\label{sec:model}
\subsection{PDE dynamics: Kuramoto-Sivashinsky equation}\label{sec:Kuramoto-Sivashinsky}
We model differential equations using a scalar time series of the Kuramoto-Sivashinsky equation ~\cite{christiansen97} under periodic boundary conditions: 
\begin{align*}
    \frac{\partial u}{\partial t} = - \frac{\partial^2 u}{\partial x^2} - \nu\frac{\partial^4 u}{\partial x^4} + \left(\frac{\partial u}{\partial x}\right)^2, \mathbb{T}\times \mathbb{R}^+
\end{align*}
where $\mathbb{T} = [0,2\pi)$, and $\nu$ corresponds to the viscosity parameter and is set as $\nu = 0.02150$. 
For creating a training time series, we use the Fourier spectral method with $32$ modes, and obtain the following form of ODEs by assuming some symmetry~\cite{christiansen97}:
\begin{align*}
&\frac{d a_k}{dt}=(k^2-\nu k^4)a_k+\frac{k}{2}\Biggl(\displaystyle\sum^{-1}_{m=k-32}a_{-m}a_{k-m}
\\&~~~~
\displaystyle\sum^{k-1}_{m=1}a_m a_{k-m}+\displaystyle\sum^{32}_{m=k+1}a_m a_{m-k}\Biggr), ~k=1,\ldots, 32.
\end{align*}

\begin{figure}
\vspace*{2mm}
  \begin{center}
\includegraphics[width=0.95\columnwidth,height=0.45\columnwidth]{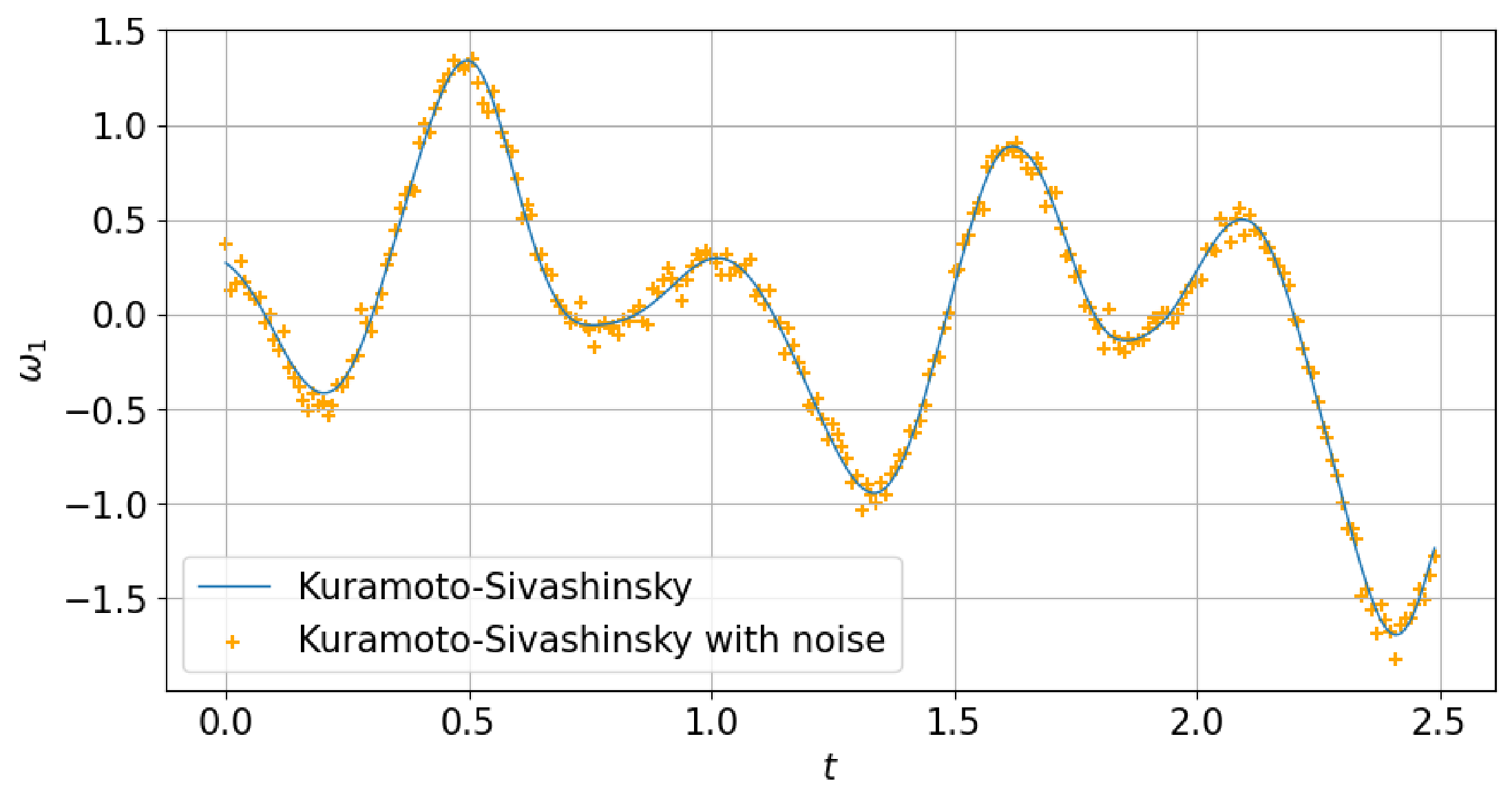}
\includegraphics[width=0.95\columnwidth,height=0.45\columnwidth]{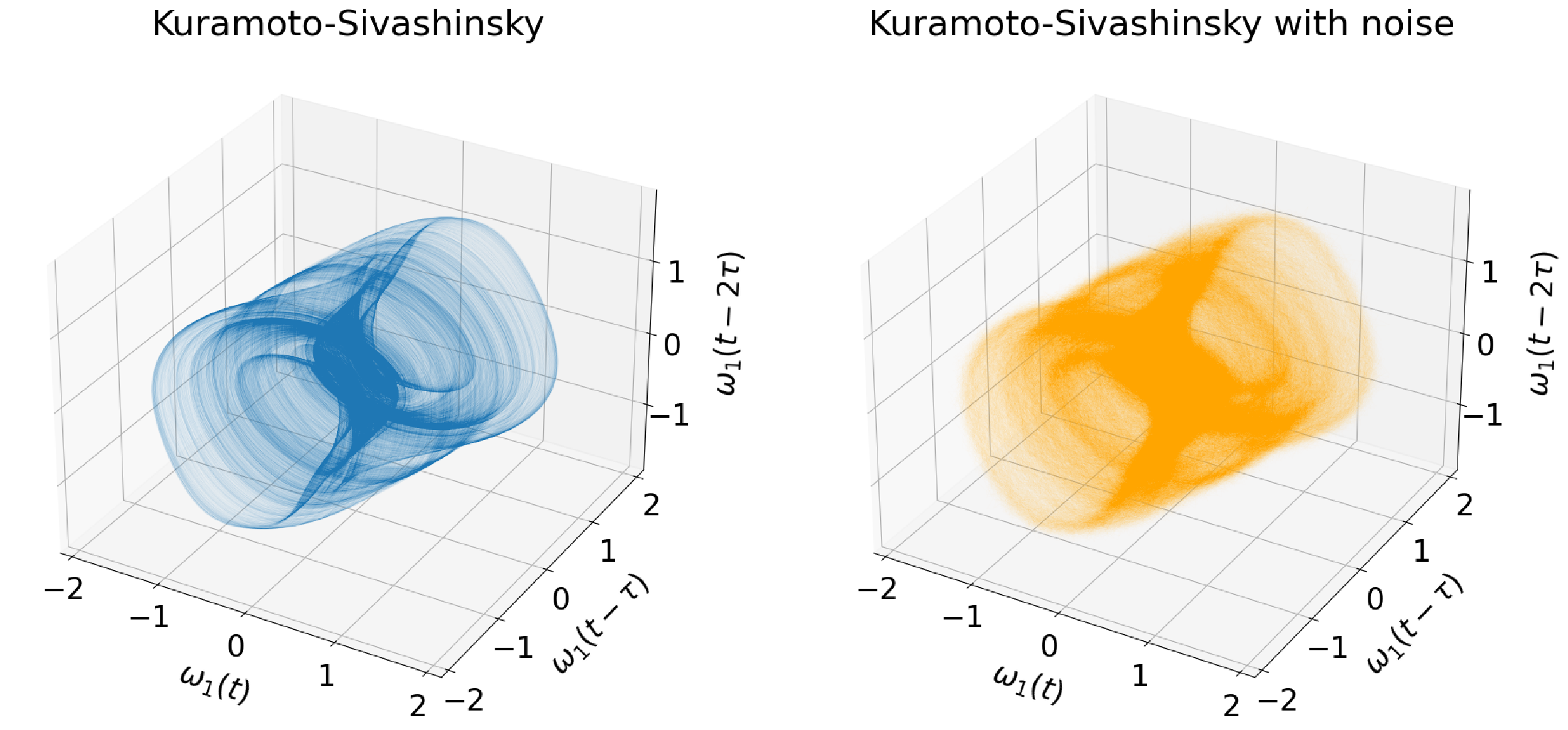}
   \end{center}
 \caption{
  {\bf Time series of $\omega_1$ of the Kuramoto-Sivasinsky equation and that with observation noise.
 The lower panels show corresponding longer trajectories in  $(\omega_1(t),\omega_1(t-\tau),\omega_1(t-2\tau))$ coordinate.}  
 The blue line represents the time series of $\omega_1$ of the Kuramoto-Sivashinsky equation and the orange points represent noised observable data.
 }
\label{fig:noised-data_KS}
\end{figure}

We assume the number of observable variables in~Eq.~\ref{eq:observable} is $I=1$, and the variable is $\omega_1 = a_1(t)$, where $a_1(t)$ is a variable of the original Kuramoto-Sivashinsky equation.
We set the model coordinate $\bm{X}$ in~Eq.~\ref{eq:odeestimation} as $\bm{X}(t) = (\omega_1(t),\omega_1(t-\tau),\ldots,\omega_1(t-4\tau))$, where $\tau =0.12$.
See Fig.~\ref{fig:noised-data_KS} for a short time series and a long trajectory in 
$(\omega_1(t),\omega_1(t-\tau),\omega_1(t-2\tau))$.

\subsection{Time-delay dynamics: Mackey-Glass equation}\label{sec:Mackey-Glass}
We deal with a delay differential equation called the 
Mackey-Glass equation~\cite{mackey77}:
\begin{align*}
    \frac{dx(t)}{dt} = 2\frac{x(t-2)}{1+x(t-2)^{9.65}} - x(t).
\end{align*}
This is also an infinite-dimensional system, but it has a finite-dimensional attractor.
We assume the number of observable variables in~Eq.~\ref{eq:observable} is $I=1$, and the variable is $\omega_1 = x(t)$, where $x(t)$ is a variable of the original Mackey-Glass equation.
We set the model coordinate as $\bm{X}$ in~Eq.~\ref{eq:odeestimation} as $\bm{X}(t) = (\omega_1(t),\omega_1(t-\tau),\ldots,\omega_1(t-6\tau))$, where $\tau =0.5$.

\subsection{Turbulence dynamics: shell model of fluid turbulence}
We model a system of differential equations using the time series of the shell model of fluid turbulence~\cite{yamada87,ohkitani89}. The system is a complex-valued differential equation, but we examine a real-valued scalar time series of an absolute value of one complex variable of the following system of differential equations: 
\begin{eqnarray*}
&&\left(\frac{du_j}{dt}+\nu k_j^2\right)u_j\nonumber\\
&&=i(c_j^{(1)}u_{j+2}^*u_{j+1}^*+c_j^{(2)}u_{j+1}^*u_{j-1}^*
+c_j^{(3)}u_{j-1}^*u_{j+1}^*)\nonumber\\
&&~~~+f \delta_{j,1}
,\quad (1\le j\le 9)
\end{eqnarray*}
where $*$ denotes the complex conjugate, $f$ is a time-independent force, 
$\nu$ is the kinematic viscosity, 
$\delta_{j,l}$ is Kronecker's delta ($l \in {\mathbb{N}}$), 
and $t$ is time. 
The real constants $c_j^{(1)},c_j^{(2)},c_j^{(3)}$ $(1\le j\le 9)$ are 
given as $c_j^{(1)}=k_j, c_j^{(2)}=-\delta k_{j-1}, c_j^{(3)}=
(\delta-1)k_{j-2}$ 
except for $c_1^{(2)}=c_1^{(3)}=c_2^{(3)}=c_{9-1}^{(1)}=c_{9}^{(1)}
=c_{9}^{(2)}=0$. 

By using a scalar time series of $w_1(t) = |u_3(t)|$
for 
$f = 0.005(1+i)$ and $\nu = 0.00251$, 
we construct a data-driven model 
by employing the six dimensional variable $\bm{X}(t) = (w_1(t), w_1(t-\tau), \ldots, w_1(t- 5\tau))$, where $\tau=18$.

\subsection{Intermittency dynamics: the  
coupled R\"ossler equation}\label{sec:Coupled-Rossler}
We model a  system of differential equations using intermittent time series of the coupled R\"ossler equation~\cite{pikovsky_1991}: 
\begin{align}
    \label{eq:coupled-rossler}
\left\{ \,
    \begin{aligned}
        & \frac{dx_1}{dt} = -y_1 -z_1 + \varepsilon (x_2-x_1) \\
        & \frac{dy_1}{dt} = x_1 + a y_1 \\
        & \frac{dz_1}{dt} = f + x_1 z_1 - cz_1 \\
        & \frac{dx_2}{dt} = - y_2 -z_2 + \varepsilon (x_1-x_2) \\
        & \frac{dy_2}{dt} = x_2 + a y_2 \\
        & \frac{dz_2}{dt} = f + x_2 z_2 - cz_2, 
    \end{aligned}
    \right.
\end{align}
where $(a, c, f, \varepsilon)=(0.15,10,0.2,0.06)$. 
This system is known for the so-called on-off intermittency, nonregular switchings between laminar and bursting states. 
Each of the states is characterized by the difference between two oscillators $(x_1,y_1,z_1)$ and $(x_2,y_2,z_2)$.
To characterize the switching, 
we focus on the difference between $x_1$ and $x_2$.
Due to the nonuniformity of the value in the difference, modeling using the difference variable is difficult.
Hence, we use two training variables $x_1$ and $x_2$ as $\omega_1$ and $\omega_2$, and see the difference. 
Then we use their delayed variables
$(X_1(t), X_2(t), X_3(t), X_4(t), X_5(t), X_6(t)) = (x_1(t), x_2(t), x_1(t-\tau), x_2(t-\tau), x_1(t-2\tau), x_2(t-2\tau))$, where $\tau =0.4$.

\subsection{PDE dynamics with noise}
\label{sec:Kuramoto-Sivashinsky-with-noise}
We model using time series data of the Kuramoto-Sivasinsky equation with
the Gaussian noise whose standard deviation is  $0.5948 \cdot 0.10$, where $0.5948$ is the standard deviation of the $\omega_1$ variable of the Kuramoto-Sivasinsky equation. 
See Fig.~\ref{fig:noised-data_KS} for a short time series and a long trajectory in  
$(\omega_1(t),\omega_1(t-\tau),\omega_1(t-2\tau))$.
By adding observation noise to the Kuramoto-Sivashinsky dynamics, 
the local structures are blurred.

\section{Results}\label{sec:results}
In this section, we assess the constructed data-driven models. 
In section~\ref{sec:Results-BasicProperties}, 
we assert that we have succeeded in predicting short-time series and that the data-driven models can reconstruct 
long trajectories
and density distribution.
We also evaluate the validity of the model by confirming the delay structure of the data-driven models.
In section~\ref{sec:Results-Basin}, we assess the stability of each model chaotic invariant set: chaotic attractor or chaotic saddle. 
We concentrate on the case when the constructed data-driven model does not have a chaotic attractor 
in section~\ref{sec:Results-stagger}. 
In section~\ref{sec:Results-Intermittency}, we assert that the data-driven model for the coupled R\"ossler equation can reconstruct intermittency.

\subsection{Basic properties of the models}\label{sec:Results-BasicProperties}
\begin{figure*}%[t]
    \centering
\includegraphics[width=2.0\columnwidth,height=1.3\columnwidth]{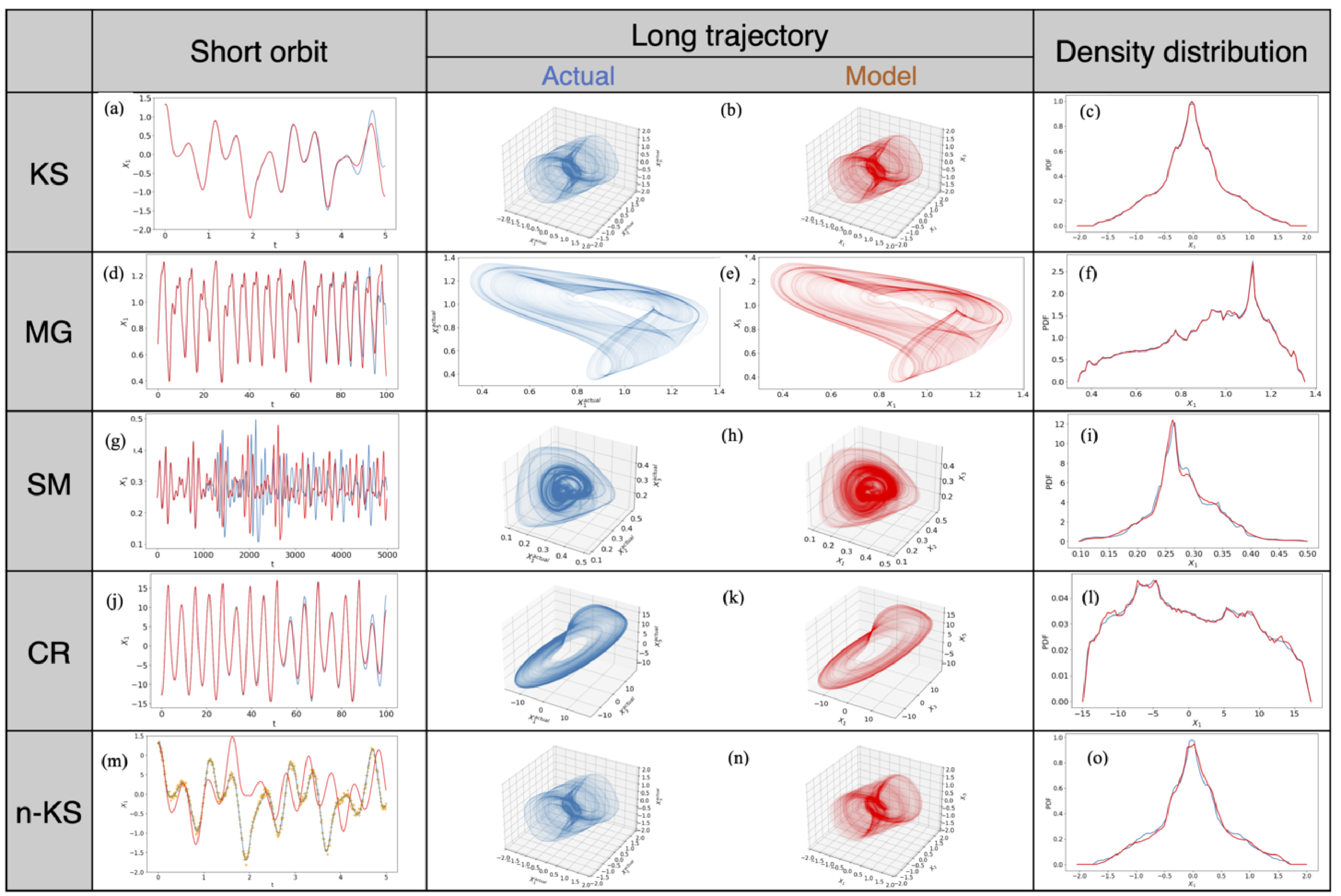}
    \caption{{\bf Basic properties of a data-driven model for each case (Kuramoto-Sivashinsky equation~(KS), Mackey-Glass equation~(MG), the shell model of fluid turbulence~(SM), the coupled R\"ossler equation~(CR), and Kuramoto-Sivashinsky equation with noise~(n-KS)).} The left panels show short trajectories of $X_1$, the center panels show 
    long trajectories,
    and the right panels show density distributions of $X_1$. In  each panel blue and red indicate the cases for the actual and the model, respectively. 
    See Appendix~\ref{sec:short-tarajectory_multi-init} for the prediction of short time orbits  from various initial conditions. 
    }
    \label{fig:all-in-one}
\end{figure*}
In this section, we evaluate the basic properties 
of each data-driven model in comparison to those of the corresponding actual dynamics.
We find that a time series inference of $X_1$ can successfully be applied for a short time. Figure~\ref{fig:all-in-one}-a,d,g,j  and m depict examples of predicted trajectories, each of which approximates the actual one for a certain amount of time. 
The growth of error in each model 
is unavoidable because of the chaotic property of the actual dynamics. 
We also confirm by calculating a long trajectory that the chaotic invariant set of each  data-driven model resembles that of the actual one. 
Figure~\ref{fig:all-in-one}-b,e,h,k  and n depict projections of long trajectories. For three data-driven models for the Kuramoto-Sivashinsky equation, the Mackey-Glass equation, and the coupled R\"ossler equation 
the obtained chaotic invariant sets are attractors,
and model trajectories are created by forward-time integration of the models.
For the case of the shell model  and the Kuramoto-Sivashinsky equation with noise,
we describe how to create an appropriate model trajectory in section~\ref{sec:Results-stagger}.
Furthermore, we demonstrate that data-driven models can reconstruct statistical quantities.
In Fig.~\ref{fig:all-in-one}-c,f,i,l and o we show a density distribution computed from a model trajectory $\{X_1(t)\}$ which coincides with that computed from the actual dynamics.
We can also see the validity of each model by confirming that short time inference from various initial conditions is successful.
See Appendix~\ref{sec:short-tarajectory_multi-init} that shows the comparison of model time series with actual ones for each of the four models. 
\begin{figure*}
    \centering   
\includegraphics[width=1.2\columnwidth,height=1.05\columnwidth]
{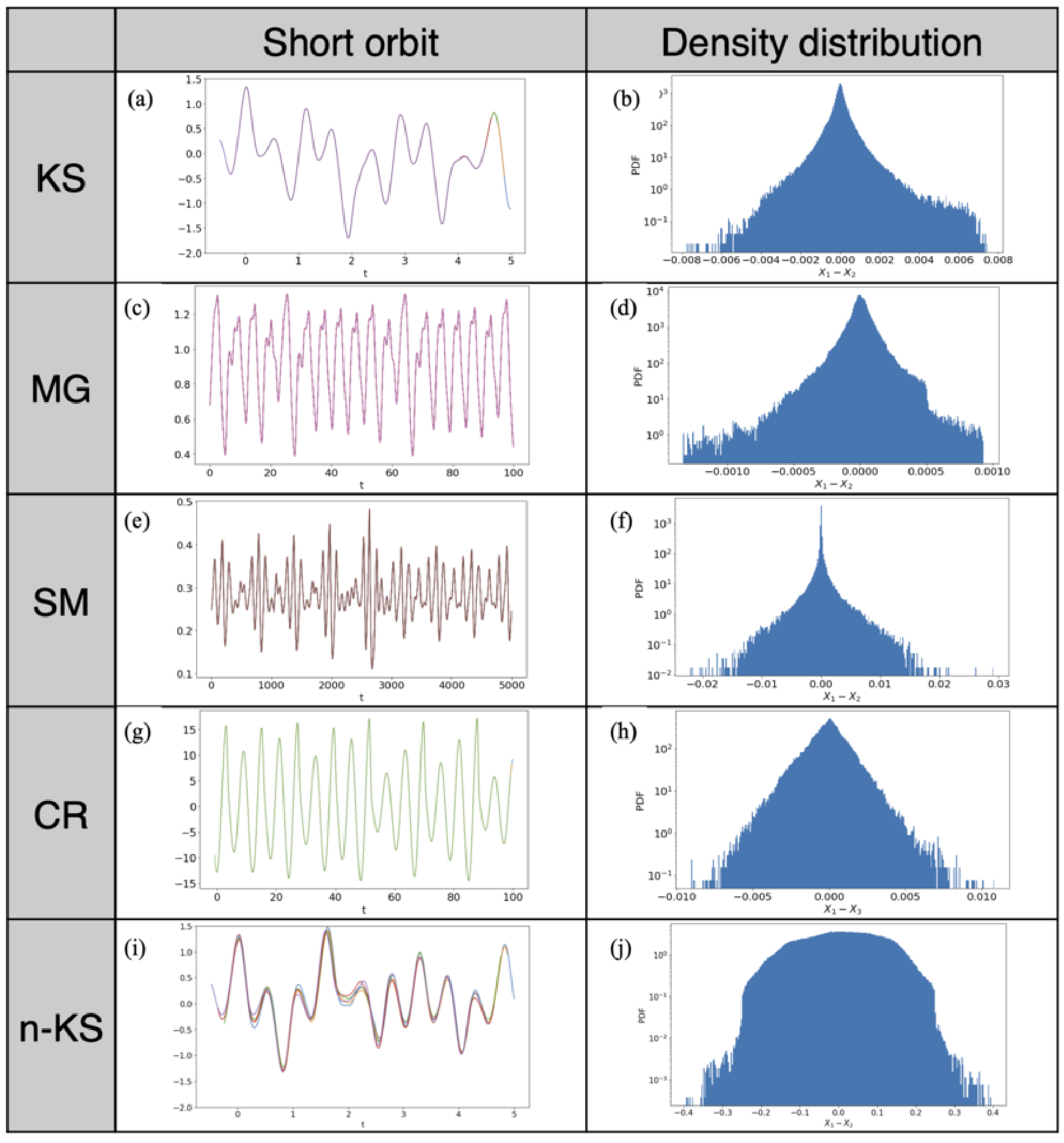}
\caption{{\bf Delay structures in each data-driven model (Kuramoto-Sivashinsky  equation~(KS), Mackey-Glass equation~(MG), the shell model of fluid turbulence~(SM), the coupled R\"ossler equation~(CR), and Kuramoto-Sivashinsky equation with noise~(n-KS)).} The left panels show short trajectories of $X_1(t)$, $X_2(t+\tau)$, $X_3(t+2\tau)$, $\cdots$, $X_{D}(t+(D-1)\tau)$ ($X_1(t)$, $X_3(t+\tau)$, $X_5(t+2\tau)$ for CR)  
    for the same time interval as that in Fig.~\ref{fig:all-in-one}~(left).
    The right panels show density distributions of $X_1(t)-X_2(t+\tau)$ ($X_1(t)-X_3(t+\tau)$ for CR) in logarithmic scale in the vertical axis, and 
  the distributions are localized around zero in comparison with the amplitudes of the fluctuations.
The results show the reproducibility of delay structures expected to satisfy for the models to be valid. 
Hence, the discrepancies in the time series in Fig.~\ref{fig:all-in-one}~(left) are considered to be due to the sensitive dependence on initial conditions of chaotic dynamics.
    }
    \label{fig:all-in-one-delay}
\end{figure*}

Recall that in our modeling we employ a delay-coordinate  variable $\bm{X}(t)$. Hence, the relation
\eqref{eq:delay_structure}
should hold for a model to be suitable.
By measuring the degree of the delay structure reconstruction, we can assess the constructed model. 
Relations shown in Fig.~\ref{fig:all-in-one-delay}-a,c,e,g and  i represent successful reconstructions of delay structure among time series  
$X_1(t)$, $X_2(t+\tau)$, $\ldots$, $X_D(t+(D-1)\tau)$ 
(for the coupled R\"ossler equation, $X_1(t)$, $X_3(t+\tau)$ and $X_5(t+2\tau)$).  
Figure~\ref{fig:all-in-one-delay}-b,d,f,h and {\allblack j} show density distributions of $X_1(t) - X_2(t+\tau)$ in each data-driven model (for the coupled R\"ossler equation, $X_1(t)-X_3(t+\tau)$).  
{\allblack
The degree of the 
reconstruction
for the Kuramoto-Sivashinsky equation with noise is lower than that of the others, due to observation noise.
}

\subsection{Stability of a chaotic invariant set}\label{sec:Results-Basin} %{\bf Basin of attraction of the model.}
Our data-driven model is constructed using trajectory points on a chaotic attractor. 
Therefore, the model cannot describe dynamics far from the attractor.
However, 
we find that the model can substantially describe dynamics outside the model attractor 
on which the trajectory points are reconstructed.

We focus on an invariant set outside the model attractor whose corresponding set does not exist in the original system. 
For many successful cases, we observe that outside the model attractor, there exists a ghost hyperbolic invariant set such as a fixed point or a periodic orbit that is $D-1$ dimensionally stable, 
and its stable manifold forms a basin boundary of the model attractor, that is, a point in the basin will be attracted to the model attractor.
In some cases, a model attractor has a global basin, and in other cases, there exist no model attractors. 

For each case of the Kuramoto-Sivashinsky  equation, the Mackey-Glass equation and the coupled R\"ossler equation a model has a chaotic attractor, 
and the attractor of a model for the coupled R\"ossler equation is a global attractor. 
 For the shell model of fluid turbulence {\allblack and the Kuramoto-Sivashinsky equation with noise}, a model attractor does not exist. 
Even though we can create a long trajectory that stays in the model  
chaotic saddle 
by using 
the Stagger-and-Step method which is described in the following section~\ref{sec:Results-stagger}.
\subsection{Stagger-and-Step method}\label{sec:Results-stagger}
As is commonly the case in the model of chaotic dynamics, a long-term trajectory 
simply created from the model does not absolutely replicate the original dynamics, 
even though the short time trajectories from multiple initial conditions 
behave adequately. 
In this section, we describe a situation and explain how to generate an appropriate long trajectory by employing the data-driven model of the shell model as an example. 

{\allblack Recall that the turbulence shell model has a chaotic attractor, and the data-driven model is constructed from a scalar time series of a trajectory on the chaotic attractor. However, the original chaotic attractor is not reconstructed as a model chaotic attractor but reconstructed as a model chaotic saddle.}

A trajectory on a chaotic saddle from almost every initial condition does not stay in the neighborhood of the chaotic saddle.
However, there exists an arbitrary long trajectory staying in the neighborhood, which we approximately create as segments of appropriate short trajectories satisfying the constraint of the delay structure explained in section~\ref{sec:delaystructure}.
% of section \ref{sec:method}
The numerical method we employ to generate such a long trajectory is  
the Stagger-and-Step method~\cite{sweet_2001c}.
We create a trajectory at least time length $T=1,000,000$ on a model chaotic saddle and confirm that the trajectory reconstructs a statistical property as well as the actual short orbit (Figs.~\ref{fig:all-in-one} and \ref{fig:all-in-one-delay}).

\allblack
In the case of the Kuramoto-Sivashinsky equation with noise, we have also succeeded in generating a model trajectory on a chaotic saddle by employing the method.
\allblack
See Appendix \ref{sec:Stagger-and-Step} for the details of the Stagger-and-Step method.
\subsection{Intermittency}\label{sec:Results-Intermittency}
\begin{figure*}[t]
\vspace*{2mm}
  \begin{center}
  \includegraphics[width=1.95\columnwidth,height=0.75\columnwidth]{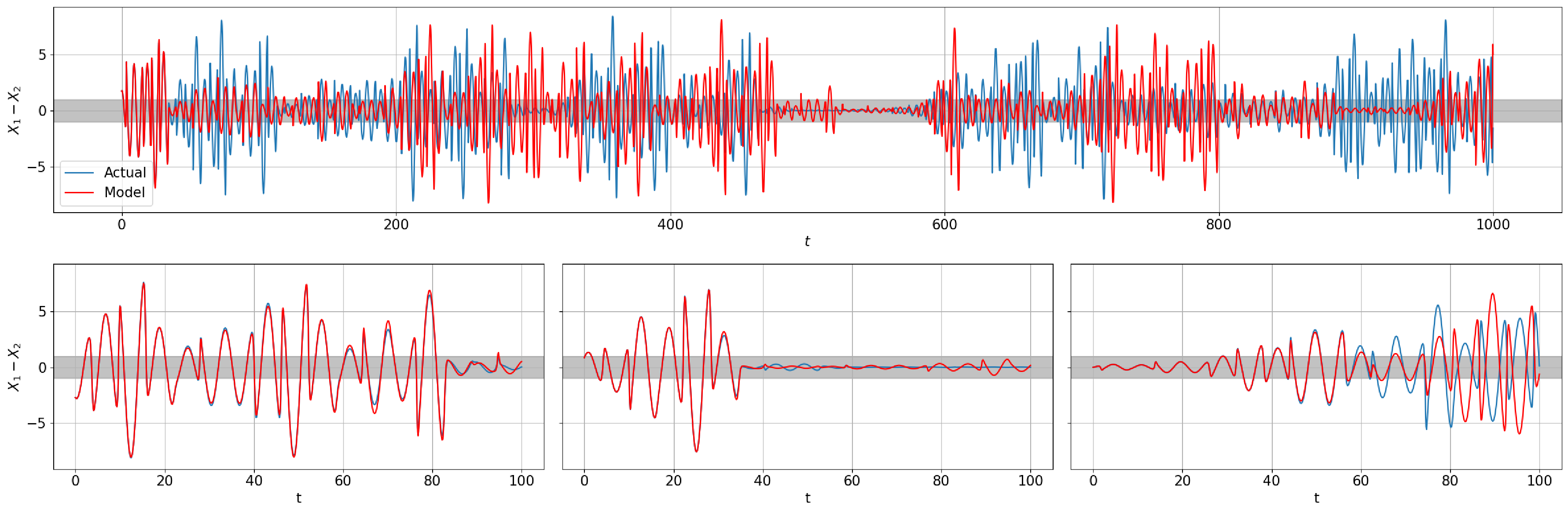}
  \end{center}
 \caption{{\bf Four time series of $X_1-X_2$ (model) from different initial conditions of a single data-driven model 
 are shown together with that of the actual dynamics (actual) of the coupled R\"ossler equation.}
 Each panel illustrates a time series showing intermittency between a laminar state with weak fluctuations $|X_1-X_2|<C~~(C=1)$ (colored in gray) and a bursting state with large fluctuations.
Each panel shows a different trajectory of a data-driven model
together with the corresponding actual trajectory. 
  Switchings between laminar and bursting states 
  can be predicted by the model.
}
\label{fig:CR_diff}
\end{figure*}
\begin{figure}
\vspace*{2mm}
  \begin{center}
  \includegraphics[width=0.95\columnwidth,height=0.65\columnwidth]
  {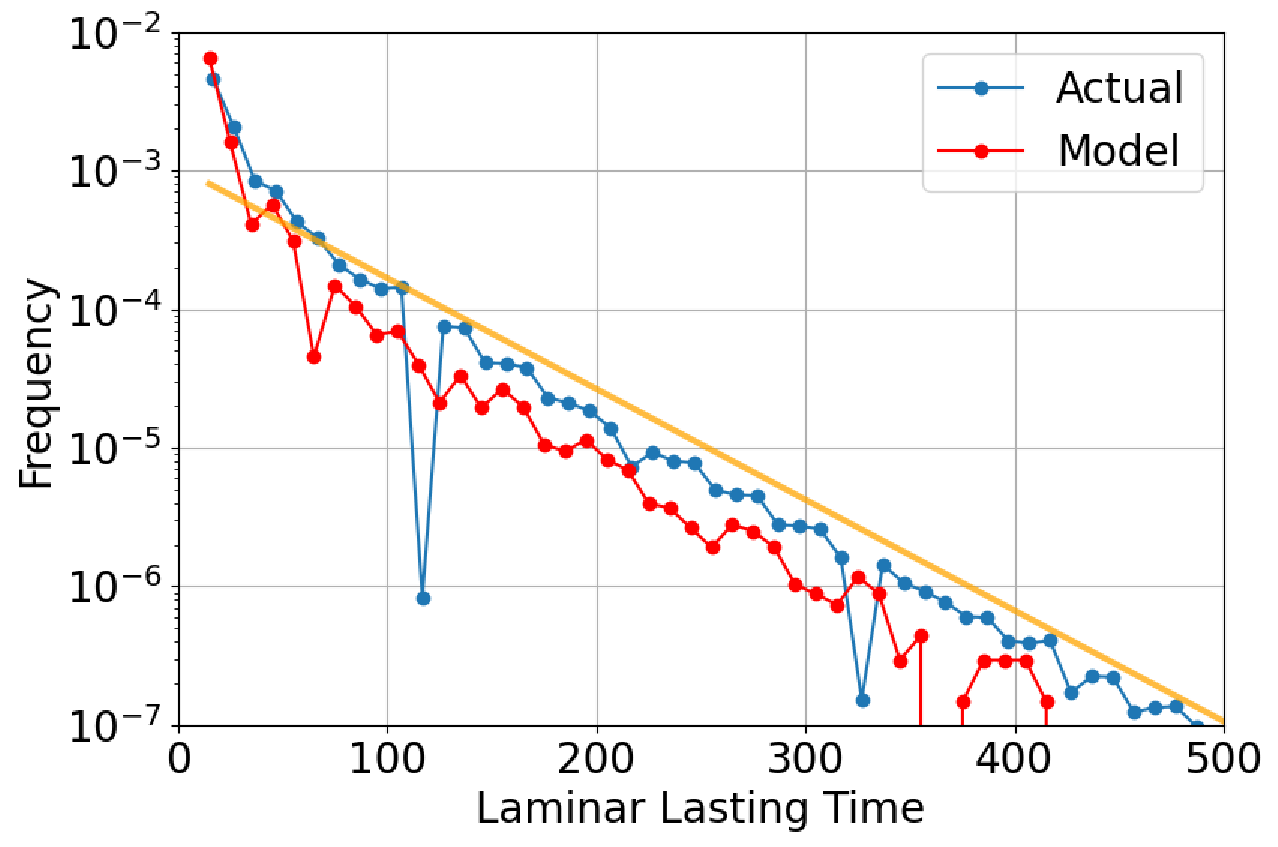}
   \end{center}
 \caption{{\bf
The distribution of laminar lasting time for the model
is shown together with that of the actual system 
in a semi-log scale for the coupled R\"ossler system.}
The Laminar lasting time distribution of $X_1-X_2$
computed from the model trajectory is displayed together with that from the actual system. 
The orange line has a slope of -0.008.
The distribution of a model is obtained from multiple numerical computations of total time length 
  $T=750,000,000$ for the actual data, and
 $T=4,330,000$ for the model.
}
\label{fig:CR_laminar_lasting_time}
\end{figure}
One of the features of the coupled R\"ossler equation~\eqref{eq:coupled-rossler} is that 
the two different oscillators can depict an intermittent behavior between laminar and bursting states through mutual interactions.
We can observe the intermittency by seeing the fluctuation of  $x_1-x_2$.

Figure~\ref{fig:CR_diff} depicts short-time trajectories of $X_1(t)-X_2(t)$ with those of the actual trajectories.
We obtain the distribution of the lasting time of the laminar state  $|x_1 - x_2|<C$ for some $C>0$ sufficiently small.
The distribution obtained by the model is depicted in Fig.~\ref{fig:CR_laminar_lasting_time} when $C=1$. 
The tail of the distribution obeys a power law, which is comparable to that from the actual coupled R\"ossler equation.
Note that it is hard to model the intermittent dynamics because of the nonuniformity when the observable variable is limited to the difference between two variables $x_1-x_2$.

\allblack 

\section{Concluding remarks}
\label{sec:conclusion}
We have succeeded in constructing various differential equations only from observable time series {\allblack of the  limited number of variables} by the RfR method proposed in our previous paper~\cite{tsutsumi22}. 
Our attempts in this study include the dynamics of 
a partial differential equation,
a delay differential equation, a shell model of fluid turbulence, and intermittent dynamics of coupled differential equations.
\allblack
The case when noise is added to the observation is also tested.
\allblack
Each data-driven model is shown to have a trajectory 
and density distribution, both of which approximate the actual ones.
The existence of the model attractor is investigated, and we create a model trajectory by applying a Stagger-and-Step method 
when an appropriate trajectory is realized not on the model attractor but on the model chaotic saddle.

\allblack In order to overcome the limitation of the number of observable variables, we employ the delay-coordinate of the observable variable(s) as a model variable, 
and it is also useful to test the validity of a model trajectory.

\allblack 
Remarkably, we can model a system of ODEs by employing the RfR method, even if 
we do not have knowledge of the system in advance.
As a practical application, we can apply our method to construct a system of ODEs from time series of macroscopic variables of complex phenomena, e.g. modeling turbulent mean flow or climate dynamics using observable macroscopic variables.
\allblack

%%%%%%%%%%%%%%%%%%%%%%%%%%%%%%%%%%%%%%%%%%%%%%%%%%%%%%%%
%%%%%%%%%%%%%%%%%%%%%%%%%%%%%%%%%%%%%%%%%%%%%%%%%%%%%%%%
\section*{acknowledgement}
YS was supported by the JSPS KAKENHI Grant No.19KK0067, 21K18584 and  23H04465, and JSPS Bilateral Open Partnership Joint Research Projects JPJSBP 120229913. 
KN was supported by the JSPS KAKENHI Grant No.22K17965 and JST PRESTO 22724051.
The computation was carried out using the JHPCN (jh220007 and jh230028) 
 and the Collaborative Research Program for Young $\cdot$ Women Scientists of ACCMS and IIMC, Kyoto University. 
%%%%%%%%%%%%%%%%%%%%%%%%%%%%%%%%%%%%%%%%%%%%%%%%%%%%%%%%%%%%%%%%%%%%%%%%%%%%%%%%%%%%%%
%%%%%%%%%%%%%%%%%%%%%%%%%%%%%%%%%%%%%%%%%%%%%%%%%%%%%%%%%%%%%%%%%%%%%%%%%%%%%%%%%%%%%%
%\bibliography{reservoir_bibliography.bib,cas-refs.bib,devarious.bib}
%
%%%%%%%%%%%%%%%%%%%%%%%%%%%%%%%%%%%%%%%%%%%%%%%%%%%%%%%%%%%%%%%%%%%%%%%%%%%%%%%%%%%%%%
%%%%%%%%%%%%%%%%%%%%%%%%%%%%%%%%%%%%%%%%%%%%%%%%%%%%%%%%%%%%%%%%%%%%%%%%%%%%%%%%%%%%%%
%%%%%%%%%%%%%%%%%%%%%%%%%%%%%%%%%%%%%%%%%%%%%%%%%%%%%%%%%%%%%%%%%%%%%%%%%%%%%%%%%%%%%%
\clearpage

\appendix

\section{Stagger-and-Step method using the delay-structures}\label{sec:Stagger-and-Step}
\begin{figure}[t]
\vspace*{2mm}
  \begin{center}
  \includegraphics[width=0.95\columnwidth,height=0.5\columnwidth]{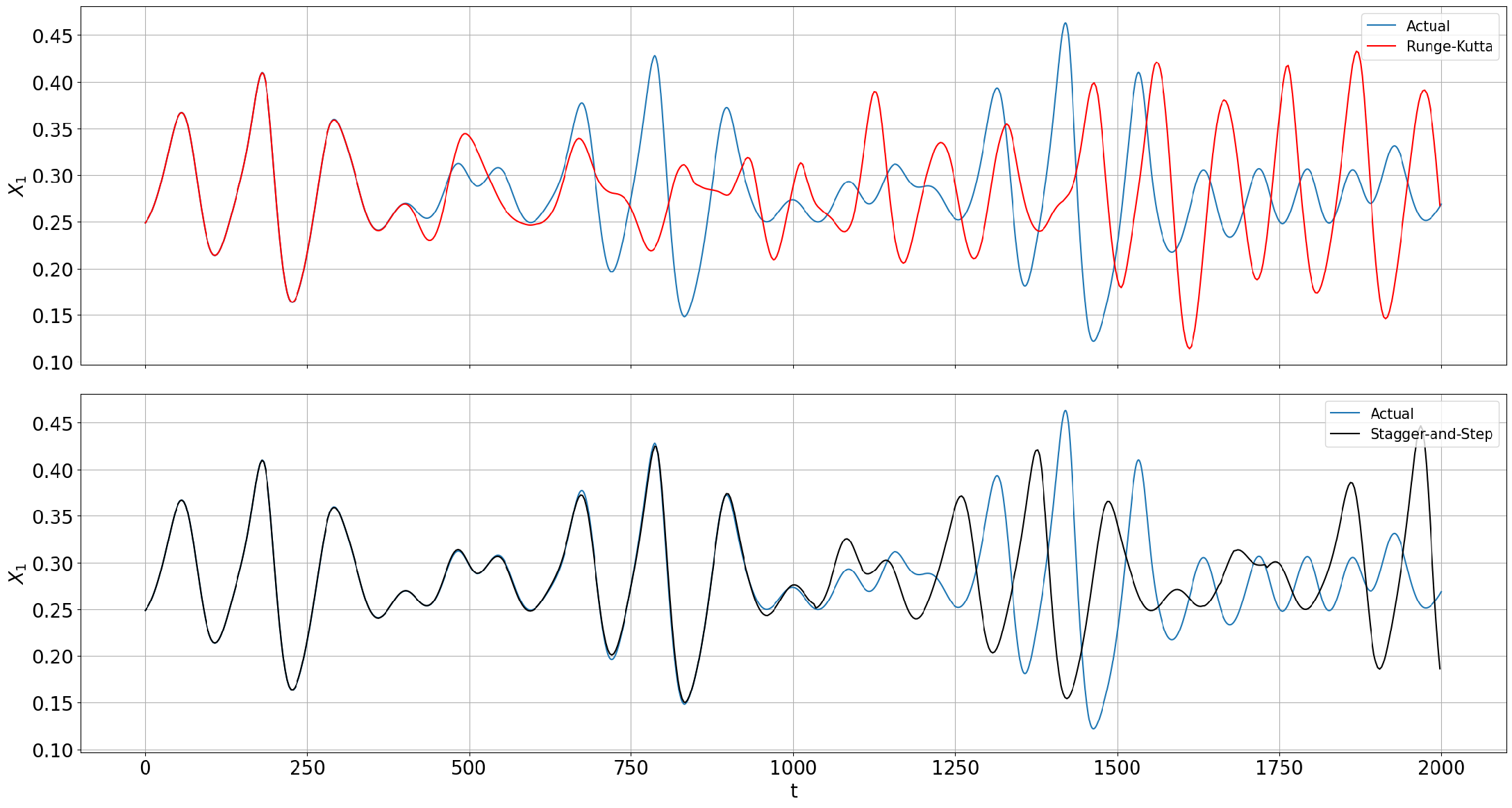}
  \end{center}
\caption{
{\bf Numerical trajectories of a data-driven model of the shell model dynamics computed from the simple calculation~(upper panel) and the Stagger-and-Step method~(lower panel).}
The blue line in each panel shows an actual trajectory. A trajectory generated using the Stagger-and-Step method is shown to approximate the actual trajectory more than times 1000, whereas a trajectory generated by simple integration of the model 
approximates the actual one for times 400.
{\allblack Note that a long trajectory created by the Stagger-and-Step method is actually a patch of short trajectories, but it seems ``smooth''. }
}
\label{fig:SM_SaS}
\end{figure}

{\allblack Stagger-and-step method~\cite{sweet_2001c} is to create a ``long trajectory'' as the patch of segments of a short trajectory.}
The method is employed to create a long trajectory on a chaotic saddle  by using the delay-structures \eqref{eq:delay_structure}.  
As seen in Fig.~\ref{fig:SM_SaS} a model trajectory created from the simple integration of the model by using the Runge-Kutta method approximates the actual one for times 400, whereas a trajectory created from the Stagger-and-Step method approximates the actual one for more than times 1000. 

We describe the Stagger-and-Step method
which we employ in Section~\ref{sec:results}.
The basic steps of the Stagger-and-Step method~\cite{sweet_2001c} are as follows.
See also  Algorithm~1.
\begin{enumerate}
    \item Calculate short numerical trajectories from nearby 100 initial conditions.
    \begin{itemize}
        \item[a] Add noise to the current point with Algorithm 2.
        \item[b] Calculate a numerical trajectory of length $T=50$.
    \end{itemize}
    \item Choose the best trajectory among 100 trajectories. Trajectories are evaluated by the maximum delay absolute error: $$E : = \displaystyle\max_{t\le 50-\tau, i=1, \ldots, D-1}|X_{i}(t) - X_{i+1}(t+\tau)|.$$
    \item Keep the former half ([0,25]) of the best trajectory as an orbit segment.
    \item Back to 1.
\end{enumerate}
To reduce the number of trials in step 1, 
\allblack
if the maximum delay absolute error $E$
is less than a certain threshold $\Delta_{\rm{threshold}}$ 
on step 1-b, we move on to step 2 even though the number of calculated numerical trajectories does not reach 100.
Furthermore, before step 1 we calculate a numerical trajectory of length $T=50$ without adding noises, and steps 1 and 2 are skipped  if $E < \Delta_{\rm{threshold}}$. 

In our computations, we slightly modify Algorithm 1 to reduce computational costs.
When the error becomes smaller with some noise, we search for a better noise by making additional noises around the original noise.
Using this optimization, an appropriate noise can be gotten with lower computational costs and a long trajectory on a chaotic saddle can be made in a plausible time.

\begin{figure}[t]
\vspace*{2mm}
  \begin{center}
  \includegraphics[width=0.9\columnwidth]{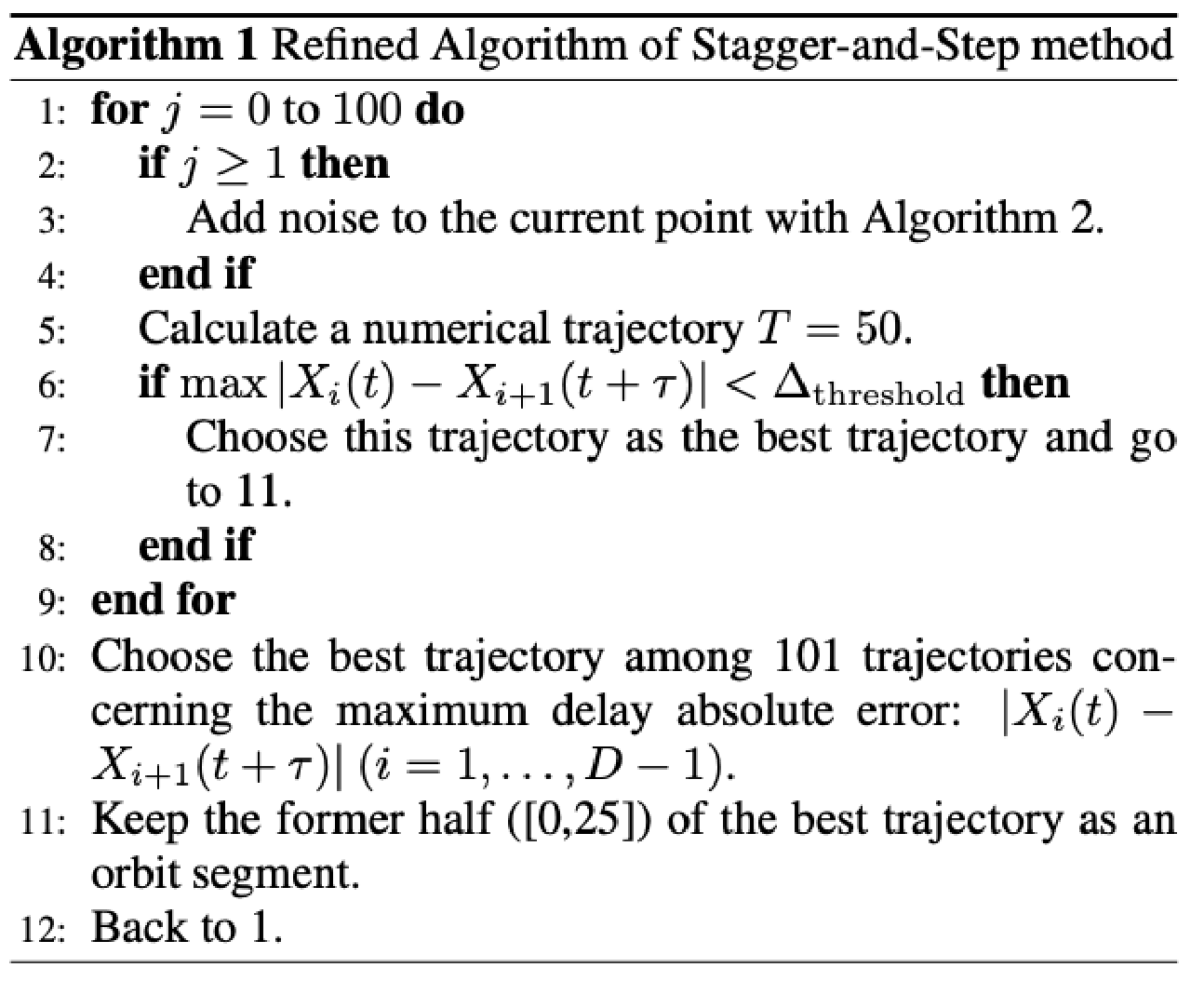}
  \end{center}
\end{figure}
\begin{figure}[t]
\vspace*{2mm}
  \begin{center}
  \includegraphics[width=0.9\columnwidth]{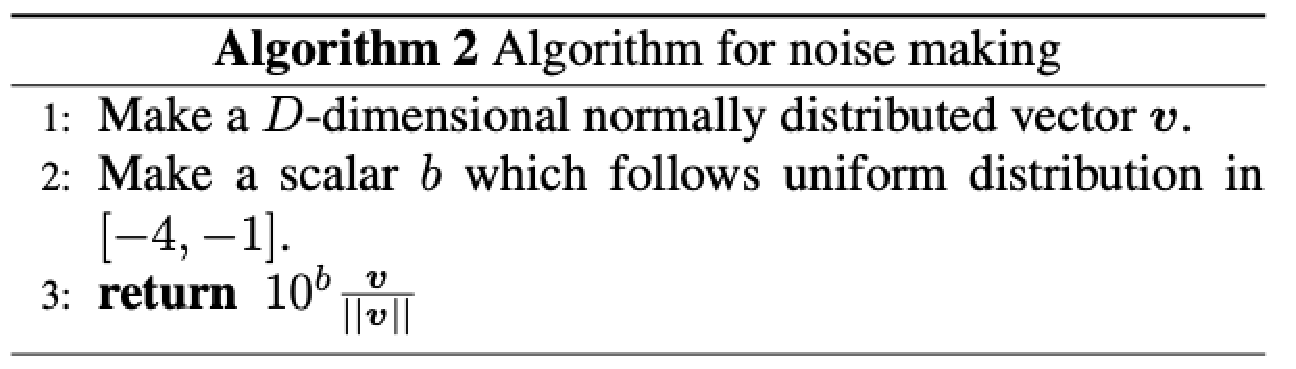}
  \end{center}
\end{figure}

% \begin{algorithm}
%  \caption{Refined Algorithm of Stagger-and-Step method}
%  \label{alg:SaS-omitting}
%  \begin{algorithmic}[1]
%   \FOR{$j = 0$ to $100$}
%   \IF{$j \ge 1$}
%   \STATE Add noise to the current point with Algorithm \ref{alg:noise-making}.
%   \ENDIF
%   \STATE Calculate a numerical trajectory $T=50$.
%   \IF{$\max |X_{i}(t) - X_{i+1}(t+\tau)| < \Delta_{\rm{threshold}}$}
%   \STATE Choose this trajectory as the best trajectory and go to 11.
%   \ENDIF
%   \ENDFOR
%   \STATE Choose the best trajectory among 101 trajectories concerning the maximum delay absolute error: $|X_{i}(t) - X_{i+1}(t+\tau)|~(i=1, \ldots, D-1)$.
%   \STATE  Keep the former half ([0,25])
%    of the best trajectory as an orbit segment.
%   \STATE Back to 1.
%  \end{algorithmic}
% \end{algorithm}

% \begin{algorithm}
%  \caption{Algorithm for noise making}
%  \label{alg:noise-making}
%  \begin{algorithmic}[1]
%   \STATE Make a $D$-dimensional normally distributed vector $\bm{v}$.
%   \STATE Make a scalar $b$ which follows uniform distribution in $[-4,-1]$.
%   \RETURN $10^b \frac{\bm{v}}{||\bm{v}||}$
%  \end{algorithmic} 
% \end{algorithm}

\clearpage
\section{Relation between auto-correlation function and the choice of delay coordinate}
\label{sec:auto-corr}
\begin{figure}[hb]
\vspace*{2mm}
  \begin{center}
  \includegraphics[width=0.95\columnwidth,height=1.2\columnwidth]{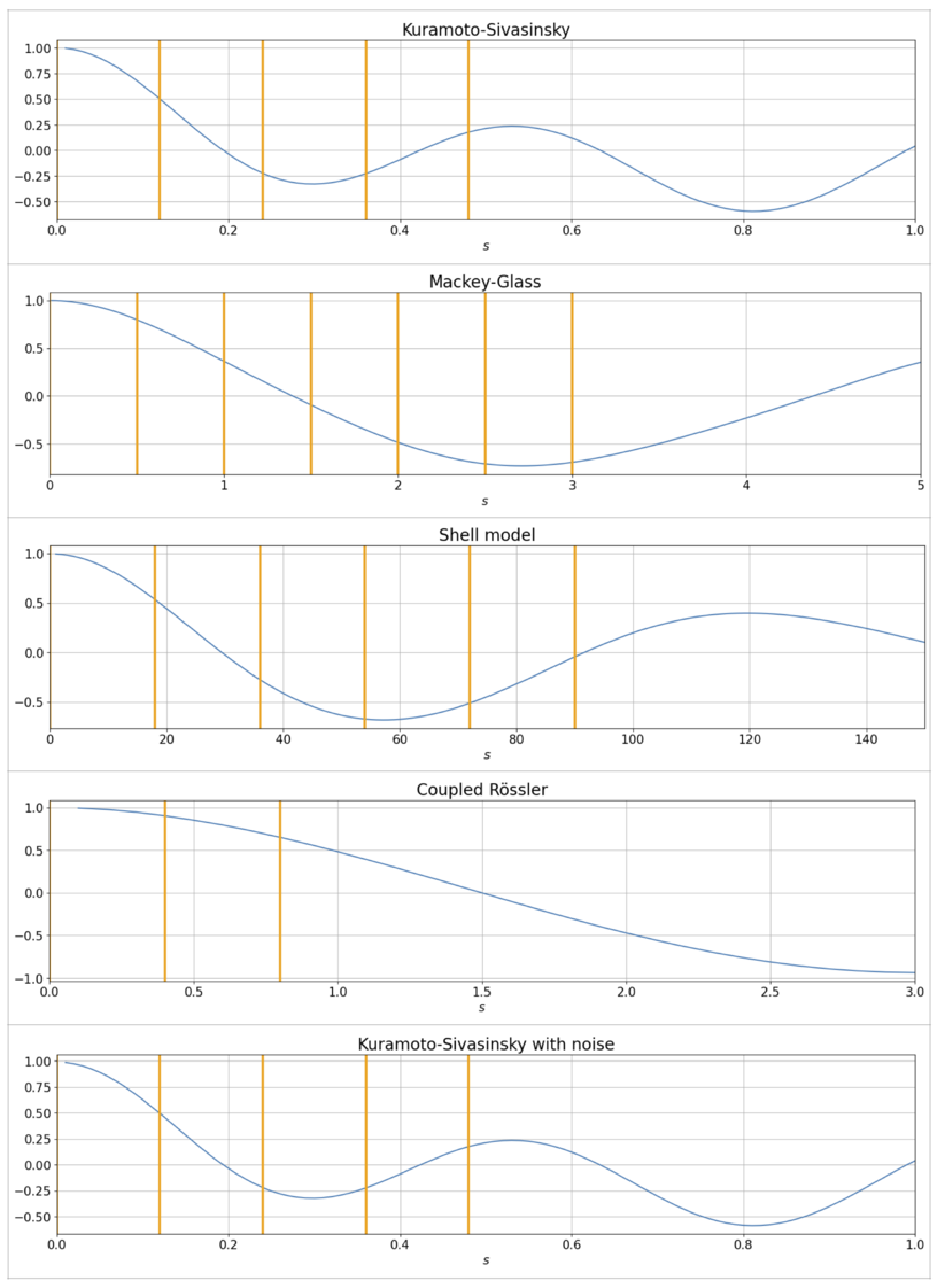}
   \end{center}
 \caption{
 {\bf Auto-correlation function of a main variable.}
 The orange lines correspond to multiples of delay time $\tau$ 
 used for the delay coordinates in Sections~\ref{sec:model} and \ref{sec:results}.
 }
\label{fig:auto-corr_all}
\end{figure}

We select the delay time $\tau$ by considering the 
decay of the auto-correlation function 
of a main variable. 
Figure~\ref{fig:auto-corr_all} depicts  auto-correlation functions of time $s$ for the observable time series data of five examples. 
The adequate choice of $\tau$ should satisfy the 
correlation coefficient between $w(t)$ and $w(t-\tau)$ to be away from both 0 and 1.
In fact, for each of the five successful examples, the correlation coefficient between 
$w(t)$ and $w(t-\tau)$ is 
KS:0.5012, MG:0.8011, SM:0.5357, CR: 0.9058, n-KS:0.4265. 

\clearpage
\allblack
\section{The estimation of time derivatives from a time series with noise}
\label{sec:relation_taylor_time-step}
To estimate time derivatives from time series data,
we employ the Taylor approximation.
As mentioned in section~\ref{sec:timederivative}.
% of Section~\ref{sec:method},
it is important to take a time step $l \Delta t$ larger when the observations include noise. 
Recall that $l \Delta t$ is a time step used to estimate time derivatives, where  $l$ is a positive integer and $\Delta t$ is a time step of observed  time series data. 

\begin{figure}[H]
\vspace*{2mm}
  \begin{center}
  \includegraphics[width=0.9\columnwidth,height=0.30\columnwidth]{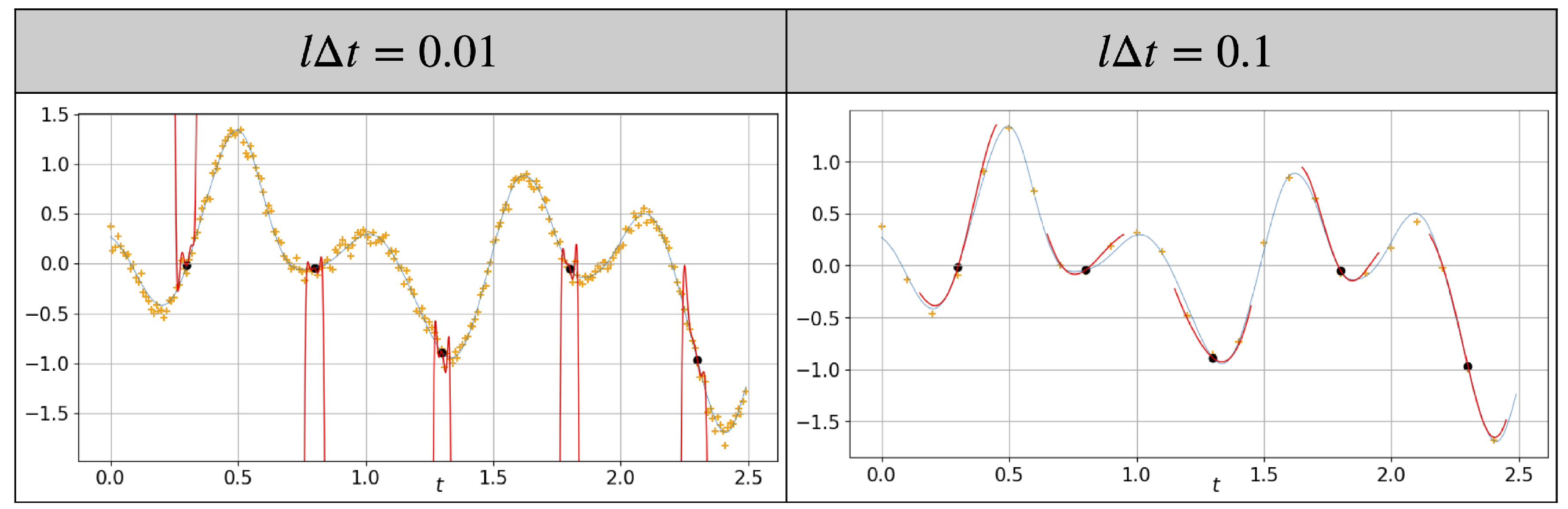}\\
\includegraphics[width=0.9\columnwidth,height=0.55\columnwidth]{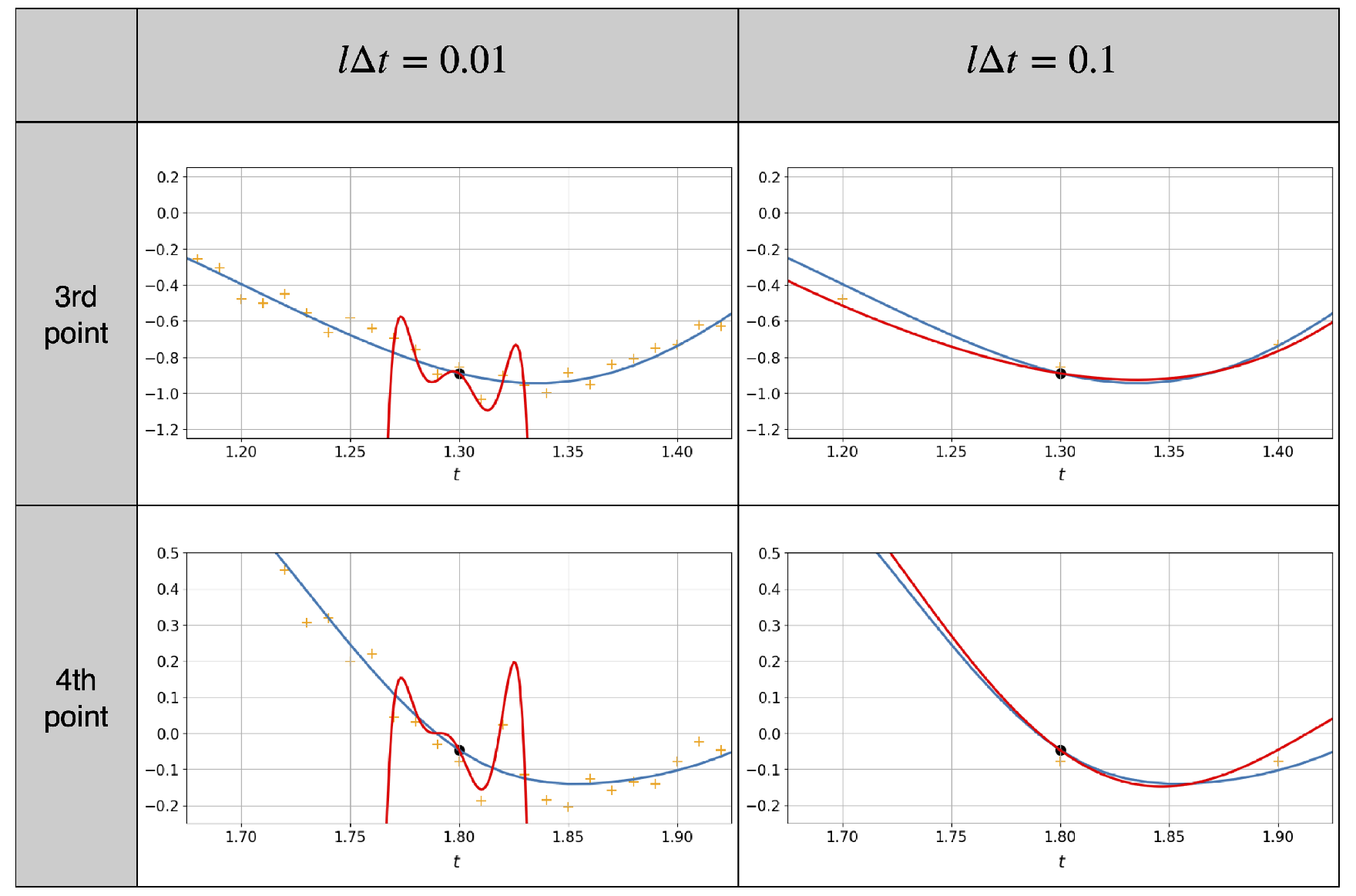}
   \end{center}
 \caption{{\allblack {\bf The sixth order Taylor approximation of a trajectory of the main variable ($\omega_1$) at some sample points for the Kuramoto-Sivashinsky equation with noise (upper) and their enlarged pictures around two sample points (lower).}
  The blue line represents the time series of the Kuramoto-Sivashinsky equation and the orange points represent noised observable data. 
  The red curve represents the Taylor approximation at each sample point 
  $( \bullet )$. 
 The left panels employ observation points at every $l\Delta t=0.01$ time step, whereas the right panels at every $l\Delta t=0.1$ time step. 
 For each of the two sample points in the lower panel, the estimation of the time derivative for the case of
 $l\Delta t=0.1$ (right) is better than that of 
 $l\Delta t=0.01$ (left).
 In order to avoid the noise effect we need to choose a large value of $l$ to estimate a time derivative at a sample point.
 }
 }
\label{fig:relation_taylor_time-step}
\end{figure}

In Fig.~\ref{fig:relation_taylor_time-step} we show the Taylor approximations at some points of the Kuramoto-Sivashinsky equation with noise.
These approximations are used to estimate time derivatives.
The estimates with a larger time step $l\Delta t$ for $l>1$ tend to obtain better approximations under the existence of observation noise. 

\begin{figure}[b]
\vspace*{2mm}
  \begin{center}
  \includegraphics[width=0.90\columnwidth,height=0.50\columnwidth]{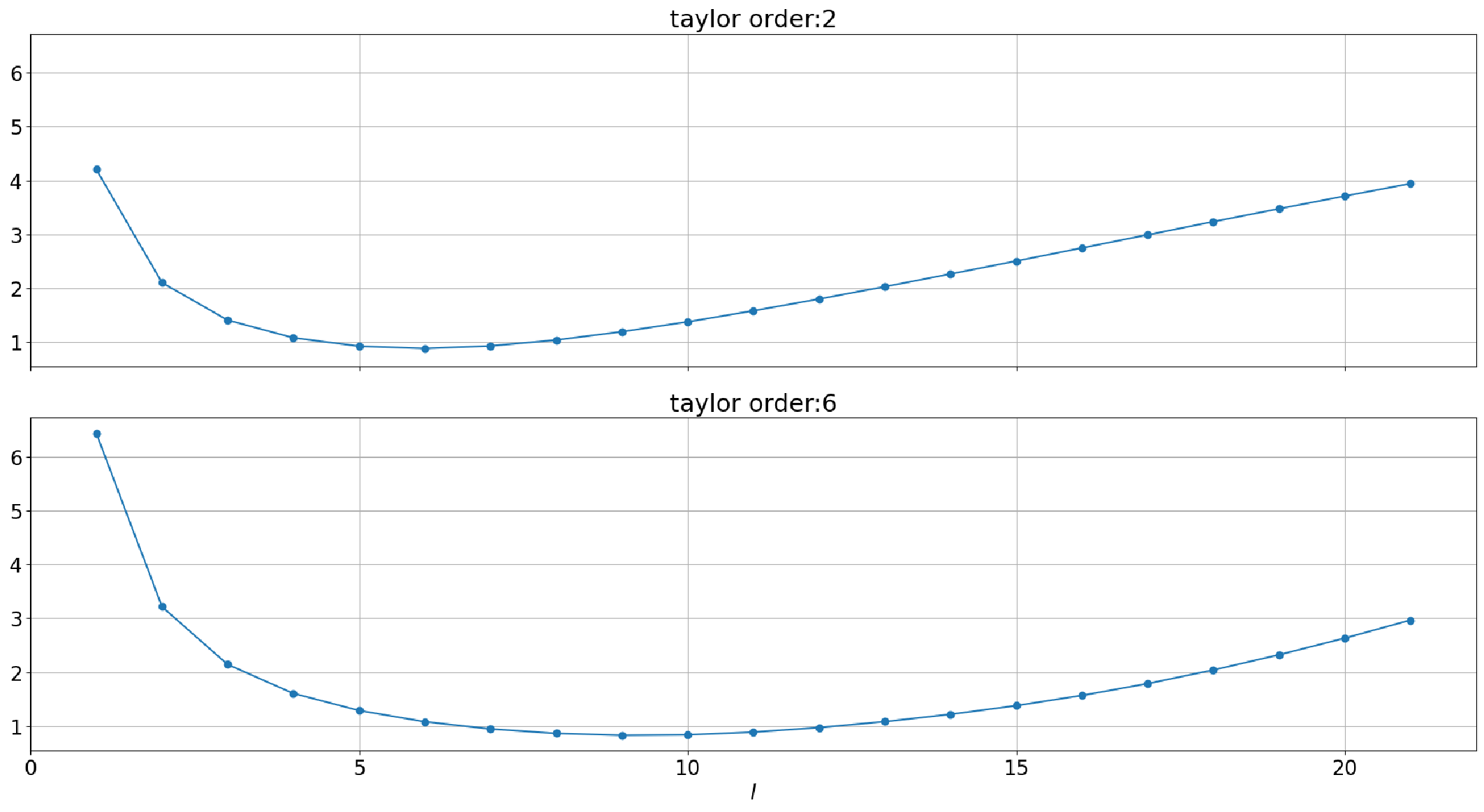}
   \end{center}
 \caption{\allblack 
 {\bf Standard deviations of the estimation errors in the time derivative with respect to $l$
used to estimate time derivatives for the two cases using the second and the sixth order Taylor approximation.}
 The two panels show the standard deviations of the estimation error of the time derivative at $10^6$ points on a trajectory. %10^6-6 
 The upper panel employs the second-order Taylor approximation, and the minimum error value $0.866$ is taken at $l = 6$.
 The lower panel employs the sixth order, the minimum error $0.833$ is taken at $l=9$.
 The standard deviation of error for the sixth-order Taylor approximation 
 tends to take a low value for a broader range of $l$.
 }
\label{fig:derivative-error_per_order}
\end{figure}

In Fig.~\ref{fig:derivative-error_per_order} we depict the standard deviations of estimation errors in the time derivatives with respect to the value of $l$. %time steps.
The minimum standard deviations of the second-order and the sixth-order approximations are 
$0.866$ and $0.833$, respectively. 
We can robustly obtain a better estimation of time derivatives for a broader range of $l$ when we apply the sixth-order approximation. 

\newpage
%\vspace{20zw}


\clearpage
\section{Short trajectories from several initial conditions }
\label{sec:short-tarajectory_multi-init}

\begin{figure}
\vspace*{2mm}
  \begin{center}
  \includegraphics[width=0.95\columnwidth,height=0.96\columnwidth]{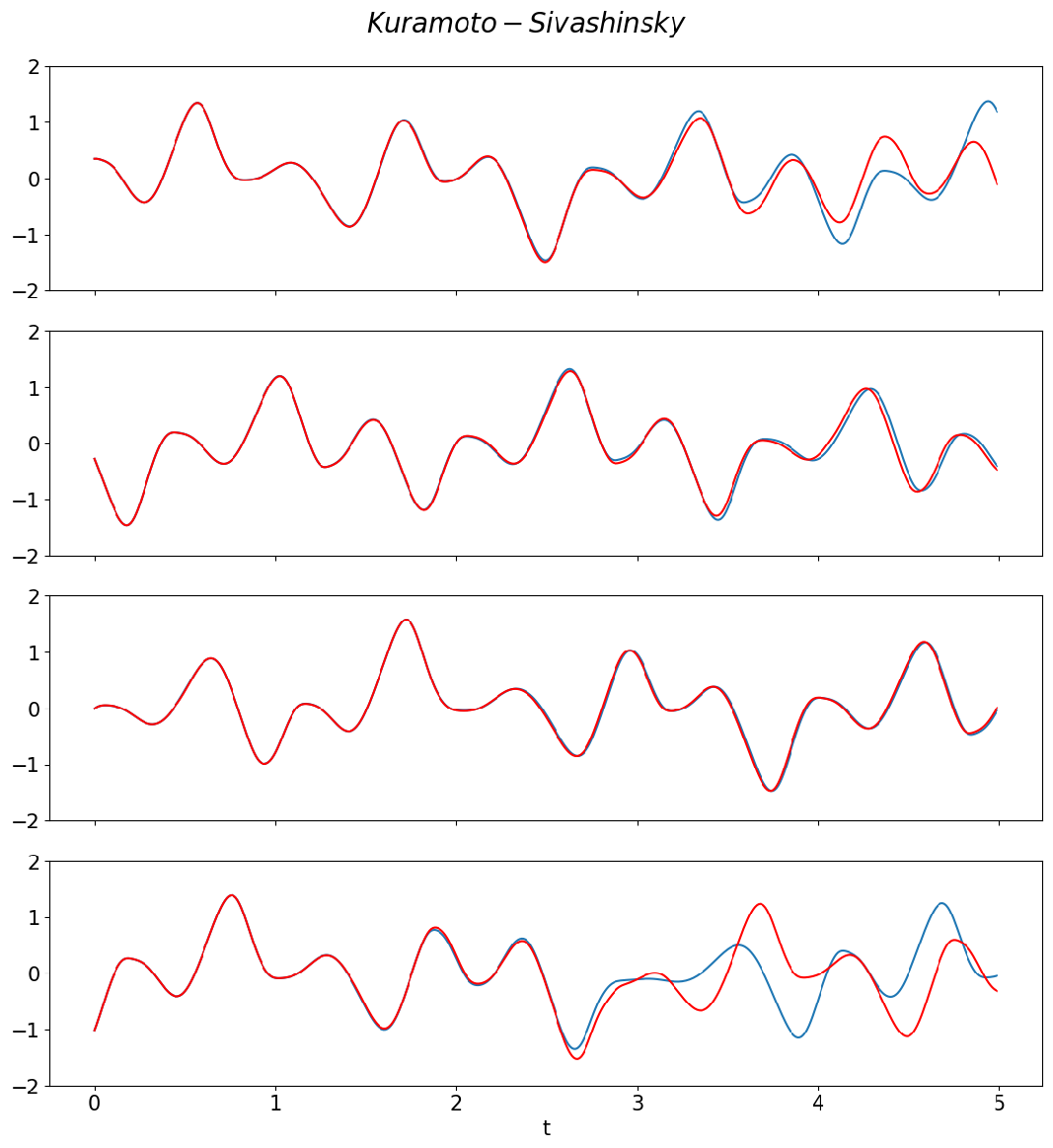}
   \end{center}
 \caption{%\allblue
 {\bf Short time trajectories of $X_1$ of a data-driven model for Kuramoto-Sivashinsky equation.}
 Blue and red indicate the cases for the actual and the model, respectively.
 Initial conditions are selected randomly from actual dynamics.
 }
\label{fig:KS_short-multi}
\end{figure}

\begin{figure}
\vspace*{2mm}
  \begin{center}
  \includegraphics[width=0.95\columnwidth,height=0.96\columnwidth]{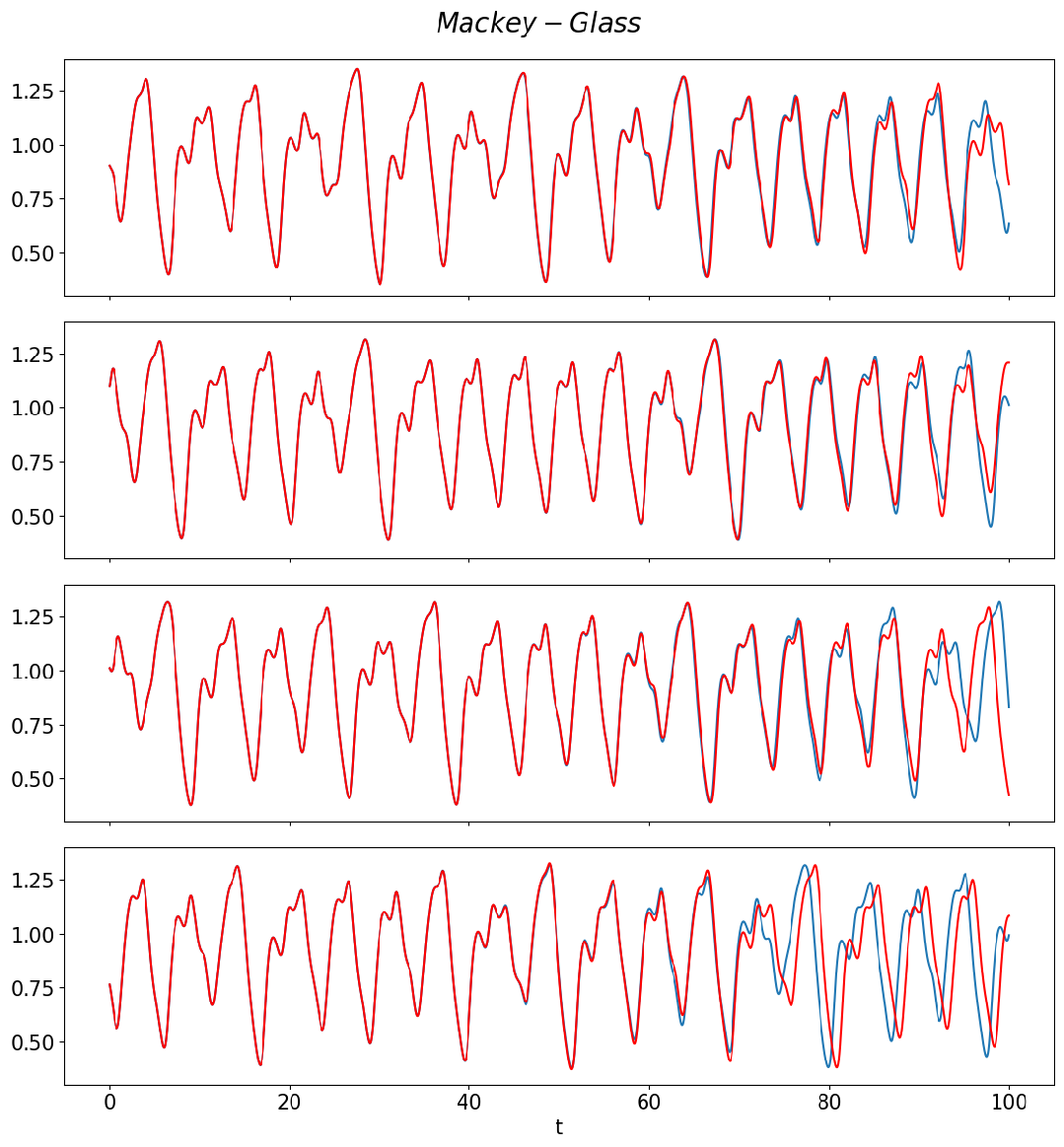}
   \end{center}
 \caption{%\allblue
 {\bf Short time trajectories of $X_1$ of a data-driven model for Mackey-Glass equation.}
 Blue and red indicate the cases for the actual and the model, respectively.
 Initial conditions are selected randomly from actual dynamics.
 }
\label{fig:MG_short-multi}
\end{figure}

\begin{figure}
\vspace*{2mm}
  \begin{center}
  \includegraphics[width=0.95\columnwidth,height=0.96\columnwidth]{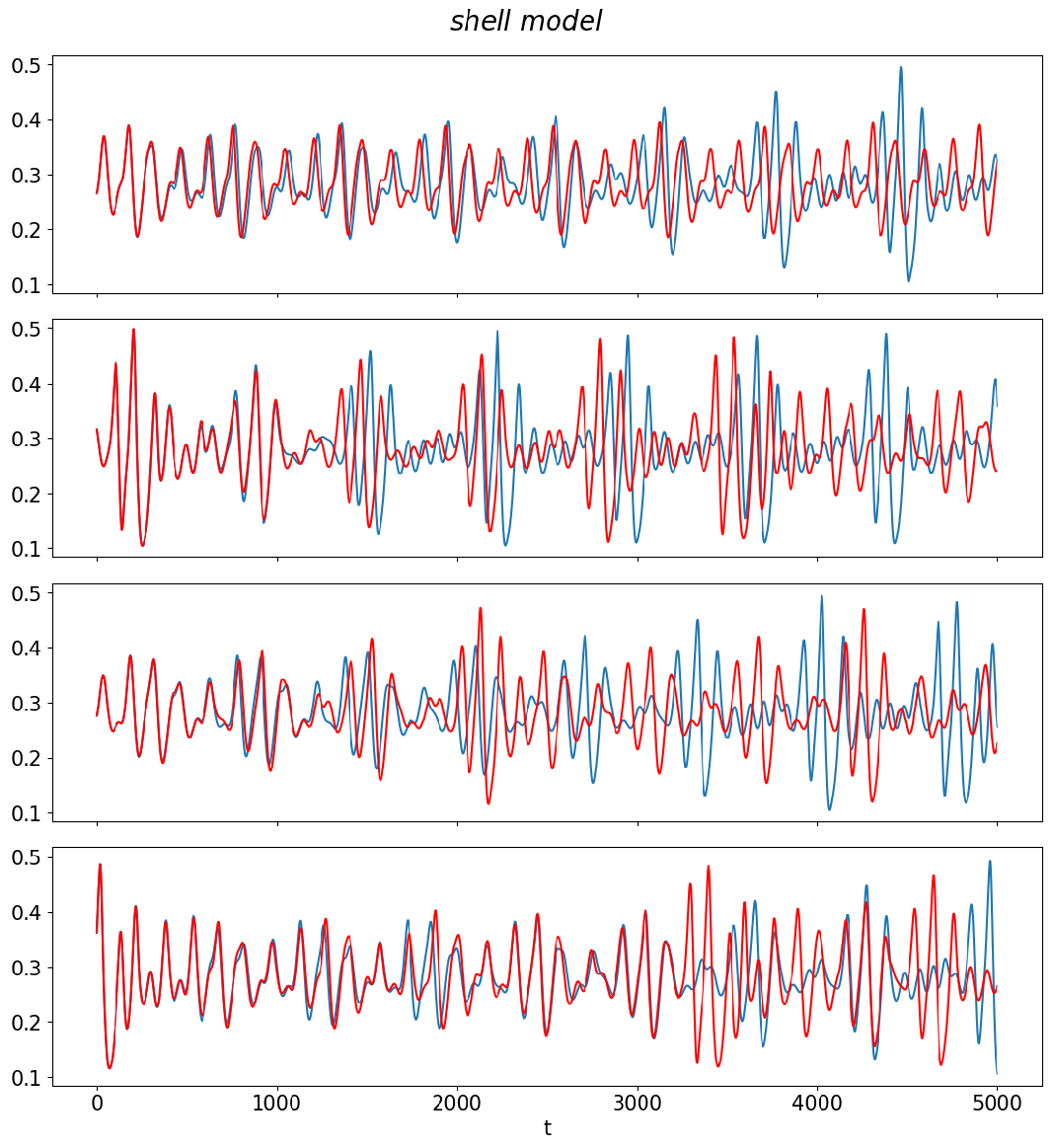}
   \end{center}
 \caption{%\allblue
 {\bf Short time trajectories of $X_1$ of a data-driven model for shell model.}
 Blue and red indicate the cases for the actual and the model, respectively.
 Initial conditions are selected randomly from actual dynamics.
 }
\label{fig:SM_short-multi}
\end{figure}

\begin{figure}
\vspace*{2mm}
  \begin{center}
  \includegraphics[width=0.95\columnwidth,height=0.96\columnwidth]{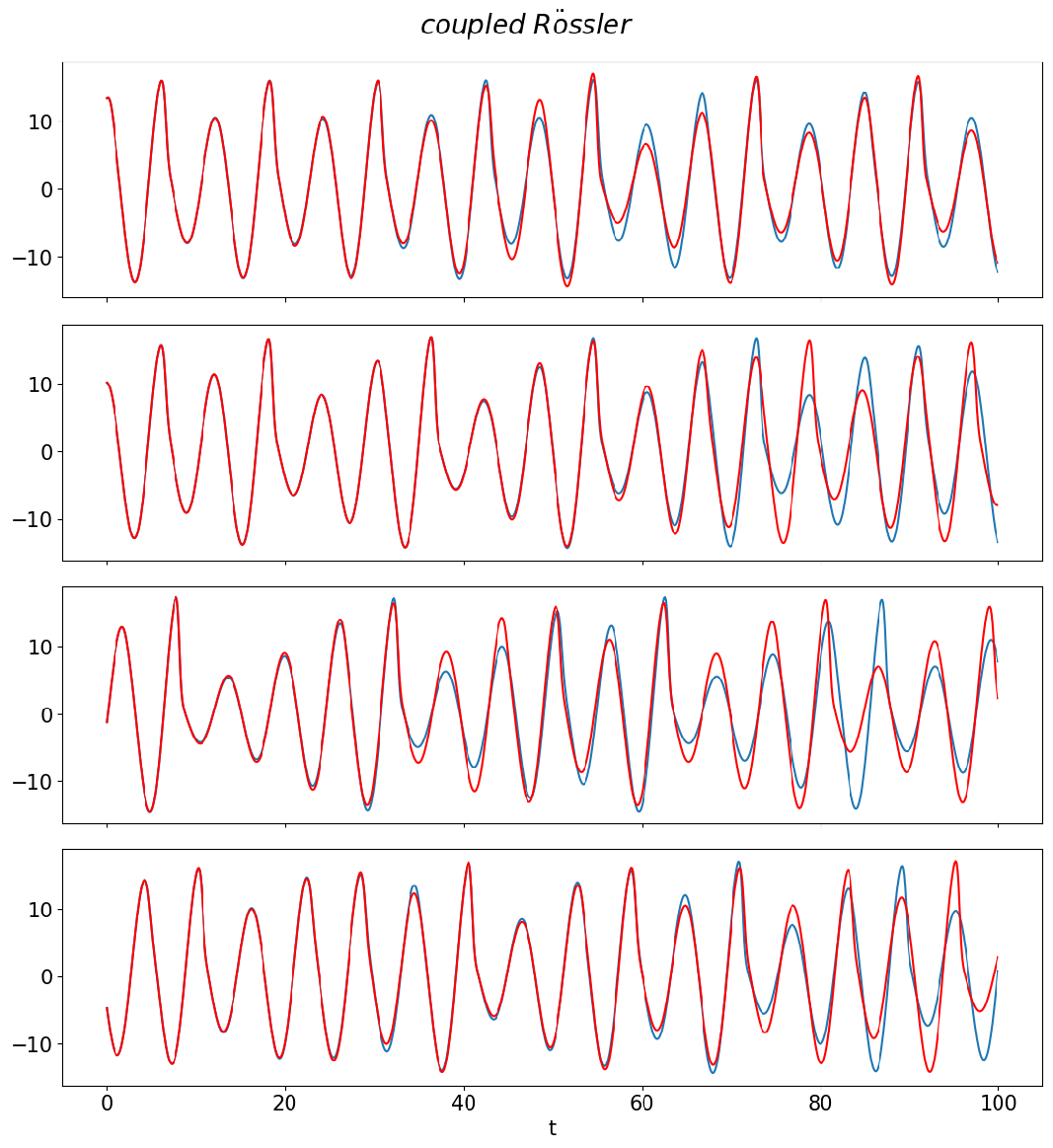}
   \end{center}
 \caption{%\allblue
 {\bf Short time trajectories of $X_1$ of a data-driven model for coupled R\"ossler equation.}
 Blue and red indicate the cases for the actual and the model, respectively.
 Initial conditions are selected randomly from actual dynamics.
 }
\label{fig:CR_short-multi}
\end{figure}

\begin{figure}
\vspace*{2mm}
  \begin{center}
  \includegraphics[width=0.95\columnwidth,height=0.96\columnwidth]{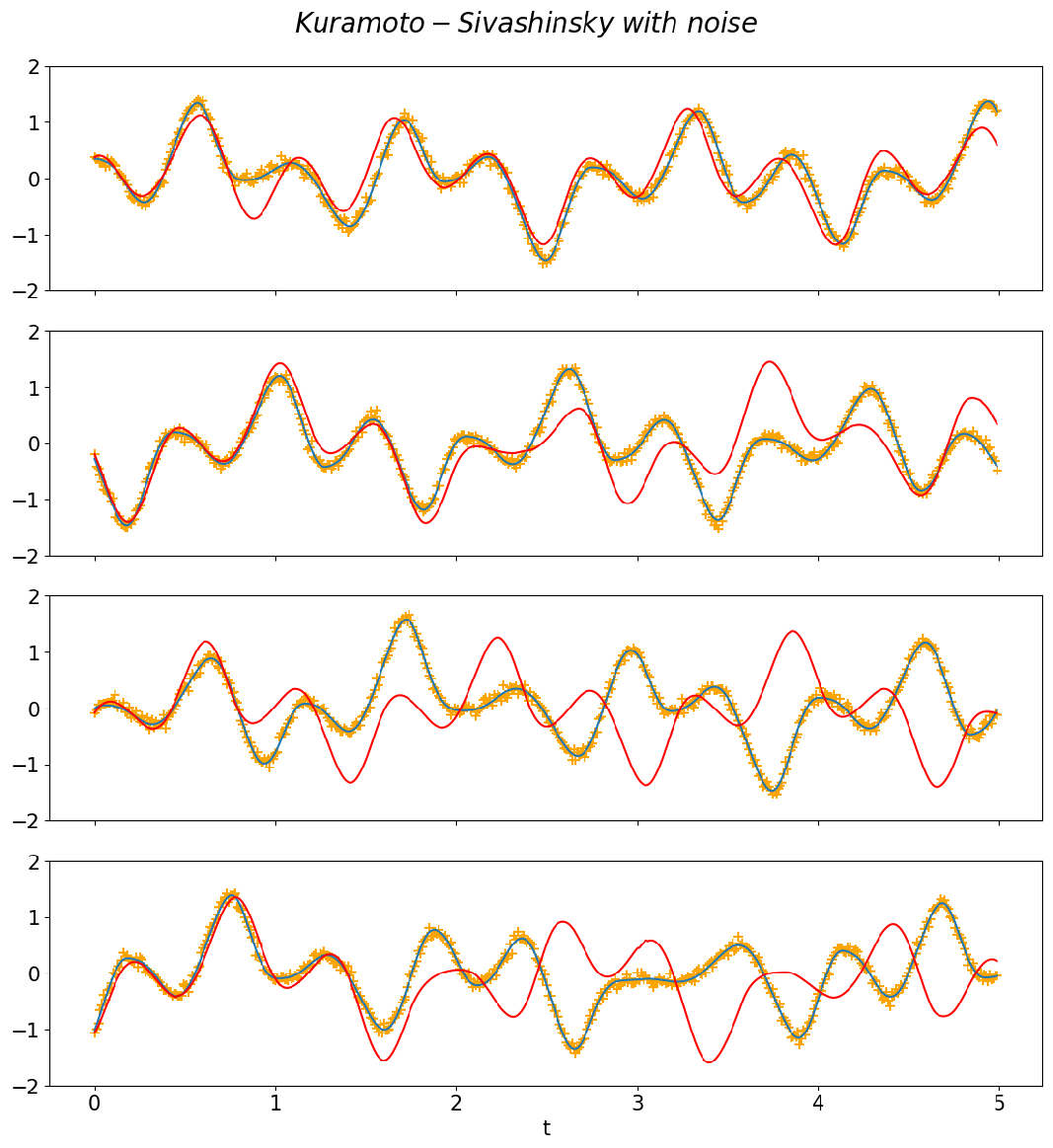}
   \end{center}
 \caption{%\allblue
 {\bf Short time trajectories of $X_1$ of a data-driven model for Kuramoto-Sivashinsky equation with noise.}
 Blue and red indicate the cases for the actual and the model, respectively and yellow points are observations with noise.
 Initial conditions are selected randomly from actual dynamics.
 }
\label{fig:n-KS_short-multi}
\end{figure}

In Section~\ref{sec:results}, we show a single short trajectory for each case in Fig.~\ref{fig:all-in-one}.
Here, we depict additional four short trajectories from randomly selected initial conditions for each of the data-driven models
in Figures~\ref{fig:KS_short-multi},\ref{fig:MG_short-multi},\ref{fig:SM_short-multi},\ref{fig:CR_short-multi}, and \ref{fig:n-KS_short-multi}.  
%to insist strongly that short inference of the data-driven models are valid.
These figures show that each model can predict short time trajectories appropriately independent of the initial conditions,
which implies the validity of each model.
% is appropriate.

\newpage

\newpage

\newpage

\newpage

\clearpage

\end{document}